\documentclass{amsart}
\usepackage{amsthm,amsmath,amsfonts,amssymb}
\usepackage{mathrsfs}
\usepackage[colorlinks=false,linktocpage=true]{hyperref}
\usepackage{tocvsec2}
\numberwithin{equation}{section}
\numberwithin{table}{section}
\font\tenscrpt=eusm10 
\font\sevenscrpt=eusm10 scaled 700
\font\fivescrpt=eusm10 scaled 500
\newfam\eusmfam
\textfont\eusmfam=\tenscrpt
\scriptfont\eusmfam=\sevenscrpt
\scriptscriptfont\eusmfam=\fivescrpt

\newtheorem{thm}{Theorem}[section]
\newtheorem{cor}{Corollary}[section]
\newtheorem{lem}{Lemma}[section]
\newtheorem{prop}{Proposition}[section]

\theoremstyle{definition}
\newtheorem{defn}{Definition}[section]

\newtheorem{rem}{Remark}[section]
\newtheorem{notn}{Notation}[section]
\newcommand{\thmref}[1]{Theorem~\ref{#1}}
\newcommand{\notnref}[1]{Notation~\ref{#1}}

\newcommand{\lemref}[1]{Lemma~\ref{#1}}
\newcommand{\coref}[1]{Corollary~\ref{#1}}
\newcommand{\propref}[1]{Proposition~\ref{#1}}
\newcommand{\defnref}[1]{Definition~\ref{#1}}
\newcommand{\remref}[1]{Remark~\ref{#1}}

\def\qed{\quad\vcenter{\hrule\hbox{\vrule height.6em\kern.6em\vrule}\hrule}}
\newenvironment{pf}{{\bigskip\textit{\newline Proof.}\quad}}{$\qed$\bigskip\newline}
\newenvironment{pf*}[1]{{\bigskip\textit{\newline#1.}\quad}}{$\qed$\bigskip\newline}

\def\ds{\displaystyle}

\def\sF{{\mathscr F}}

\def\sFt{{\mathscr F}_t}

\def\OFP{(\Omega,\sF,\P)}
\def\OFFtP{(\Omega,\sF,\{\sFt\},\P)}

\def\R{\mathbb R}
\def\E{\mathbb E}
\def\P{{\mathbb P}}
\def\Pt{\widetilde\P}
\def\Et{\widetilde\E}

\def\P{{\mathbb P}}
\def\Pt{\tilde{\mathbb P}}

\def\E{{\mathbb E}}

\def\R{\mathbb R}

\def\S{\mathbb S}
\def\T{\mathbb T}

\def\Rp{{\R}_+}

\def\lbl#1{\label{#1}}

\def\lqv{\left<}
\def\rqv{\right>}

\def\lab{\left|}
\def\rab{\right|}
\def\lgab{\Big{|}}
\def\lpa{\left(}
\def\rpa{\right)}
\def\lbk{\left[}
\def\rbk{\right]}
\def\lbr{\left\{}
\def\rbr{\right\}}
\def\bdf{\begin{defn}}
\def\edf{\end{defn}}
\def\bcr{\begin{cor}}
\def\ecr{\end{cor}}
\def\bnt{\begin{notn}}
\def\ent{\end{notn}}
\def\brm{\begin{rem}}
\def\erm{\end{rem}}
\def\blm{\begin{lem}}
\def\elm{\end{lem}}
\def\bpf{\begin{pf}}
\def\bpfs{\begin{pf*}}
\def\epf{\end{pf}}
\def\epfs{\end{pf*}}

\def\beq{\begin{equation}}
\def\beqs{\begin{equation*}}
\def\eeq{\end{equation}}
\def\eeqs{\end{equation*}}
\def\bsp{\begin{split}}
\def\esp{\end{split}}
\def\bc{\begin{cases}}
\def\ec{\end{cases}}
\def\bt{\begin{tabular}}
\def\et{\end{tabular}}

\def\bthm{\begin{thm}}
\def\ethm{\end{thm}}
\def\bpr{\begin{prop}}
\def\epr{\end{prop}}
\def\babs{\begin{abstract}}
\def\eabs{\end{abstract}}
\def\lbl{\label}
\def\bitm{\begin{itemize}}
\def\eitm{\end{itemize}}
\def\ben{\begin{enumerate}}

\def\rencomrom{\renewcommand{\labelenumi}{(\roman{enumi})}}
\def\een{\end{enumerate}}
\def\Rp{{\mathbb R}_{+}}

\def\sF{{\mathscr F}}
\def\sDot{{\mathscr D}_{1,2}}
\def\sFs{{\mathscr F}_s}
\def\sFt{{\mathscr F}_t}
\def\sFtp{{\mathscr F}_{t^{+}}}
\def\sFto{{\mathscr F}_{t_{1}}}
\def\sFtt{{\mathscr F}_{t_{2}}}
\def\sFtk{{\mathscr F}_{t_{k}}}
\def\sFtn{{\mathscr F}_{t_{n}}}
\def\sFtnpo{{\mathscr F}_{t_{n+1}}}
\def\sFT{{\mathscr F}_T}
\def\sFTt{\widetilde{\mathscr F}_T}
\def\OFTP{(\Omega,\sFT,\P)}
\def\OFFtP{(\Omega,\sF,\{\sFt\},\P)}
\def\OFTFtPt{(\Omega,\sFT,\{\sFt\},\Pt)}

\def\WrT{W_{|_{[0,T]}}}
\def\WrTn{W_{|_{[0,T_{0}]}}}
\def\sPtworT{\mathscr{P}_{2}\lpa\WrT\rpa}

\def\sPtwosrT{\mathscr{P}_{2}^{\mbox{pr}}\lpa\WrT\rpa}

\def\sPtwolocrT{\mathscr{P}_{2}^{\mbox{loc}}(\WrT)}

\def\sPtwoslocrT{\mathscr{P}_{2}^{\mbox{pr},\mbox{loc}}(\WrT)}
\def\sPtwoslocrTn{\mathscr{P}_{2}^{\mbox{pr},\mbox{loc}}(\WrTn)}
\def\DW{\mathbb {D}_{W}} 
 
\def\DWtld{\mathbb {D}_{\widetilde{W}}} 
\def\DWs{\mathbb {D}_{W_{s}}} 
\def\DWt{\mathbb {D}_{W_{t}}}
\def\DWEFsF{\mathbb{D}_{W}\mathbb{E}\left[ F|\mathscr{F}\right]}
\def\DWtEFsFt{\mathbb{D}_{W_{t}}\mathbb{E}\left[ F|\mathscr{F}_{t}\right]}
\def\DWEFsFaie{\mathbb{D}_{W}^{\sc{aie}}\mathbb{E}\left[ F|\mathscr{F}\right]}

\def\lnrm{\left\Vert}
\def\rnrm{\right\Vert}
\def\Wt{\widetilde{W}}

\def\vfi{\varphi}
\def\tone{{t_{1}}}
\def\tt{{t_{2}}}
\def\tk{{t_{k}}}
\def\tn{{t_{n}}}
\def\tnpo{{t_{n+1}}}

\def\DWtone{\mathbb {D}_{W_{\tone}}}
\def\DWttone{\mathbb {D}_{\widetilde{W}_{\tone}}}

\def\DWtt{\mathbb {D}_{W_{\tt}}}
\def\DWttt{\mathbb {D}_{\widetilde{W}_{\tt}}}
\def\DWttn{\mathbb {D}_{\widetilde{W}_{\tn}}}
\def\DWtk{\mathbb {D}_{W_{\tk}}}
\def\DWtn{\mathbb {D}_{W_{\tn}}}
\def\DWtnpo{\mathbb {D}_{W_{\tnpo}}}
\def\ft{\hat{f}}

\def\ds{\displaystyle}

\def\dsty{\displaystyle}
\begin{document}
\title[Applications of the quadratic covariation differentiation theory]{Applications of the quadratic covariation differentiation theory: variants of the Clark-Ocone and Stroock's formulas}
\author{Hassan Allouba}
\thanks{Hassan Allouba is the corresponding author}
\email{allouba@math.kent.edu}
\address{Department of Mathematical Sciences, Kent State University, Kent, Ohio 44242}
\date{3/22/2011}
\author{Ramiro Fontes}
\babs
In a 2006 article (\cite{A1}), Allouba gave his quadratic covariation differentiation theory for It\^o's integral calculus.   In it,  he defined the derivative of a semimartingale with respect to a Brownian motion as the time derivative of their quadratic covariation and a generalization thereof.  He then obtained a systematic pathwise stochastic differentiation theory that comes  complete with a fundamental theorem of stochastic calculus relating this derivative to It\^o's integral, a differential stochastic chain rule, a differential stochastic mean value theorem, and other differentiation rules.  In this current article, we  use this differentiation theory in \cite{A1} to obtain variants of the celebrated Clark-Ocone and Stroock representation formulas, with and without change of measure.  We prove our variants of the Clark-Ocone formula under $L^{2}$-type conditions on the random variable but with no $L^{p}$ conditions on the derivative.  We do not use Malliavin calculus, weak distributional or Radon-Nikodym type derivatives, or the significant extra machinery of the Hida-Malliavin calculus.  Moreover, unlike with Malliavin or Hida-Malliavin  calculi, the form of our variant of the Clark-Ocone formula under change of measure is as simple as it is under no change of measure, and without requiring any further differentiability conditions on the Girsanov transform integrand beyond the standard Novikov condition.  This is a consequence of the invariance under change of measure of the first author's derivative in \cite{A1}.  The formulations and proofs are simple and natural applications of the differentiation theory in \cite{A1} and standard It\^o integral calculus. Iterating our variants of the Clark-Ocone formula, we obtain variants of Stroock's formula.     We illustrate the applicability of these formulas and the theory in \cite{A1} by easily, and without Hida-Malliavin methods, obtaining the representation of the Brownian indicator $F=\mathbb{I}_{[K,\infty)}\left( W_{T}\right)$, which is not standard Malliavin differentiable, and by applying them to digital options in finance.   We then identify the chaos expansion of the Brownian indicator.  The first author further extends and applies his differentiation theory in forthcoming articles and obtains a general stochastic calculus for a large class of processes with different orders and types of variations, including many that fall outside the classical Gaussian, Markovian, or semimartingale classes. \eabs
\maketitle
\newpage
 \tableofcontents
\section{Introduction and statement of results}
In \cite{A1} Allouba gave his quadratic covariation pathwise stochastic differentiation theory of semimartingales with respect to Brownian motion (BM).  His idea starts by defining the strong stochastic derivative $\mathbb{D}_{W_{t}}S_{t}={dS_{t}}/{dW_{t}}$ of the ``temporally-rough'' continuous semimartingale $S$ with respect to the  ``comparably temporally-rough'' Brownian motion  $W$ at time $t$ in terms of the derivative ${d\left\langle S,W\right\rangle _{t}}/{d\left\langle W\right\rangle _{t}}$ of the ``temporally-smooth'' quadratic covariation of $S$ and $W$, $\left\langle S,W\right\rangle$, with respect to the  ``comparably temporally-smooth''  quadratic variation of $W$, $\left\langle W\right\rangle$, at $t$:
 \begin{equation}\label{SAD}
\mathbb{D}_{W_{t}}S_{t}(\omega):=\frac{d\left\langle S,W\right\rangle _{t}\left(
\omega \right) }{d\left\langle W\right\rangle _{t}\left( \omega \right) }=%
\frac{d\left\langle S,W\right\rangle _{t}\left( \omega \right) }{dt},
\end{equation}%
almost surely (see Definition 1.1 equation (2) in \cite{A1} and the quadratic covariation \defnref{covXY} below).   He then develops in \cite{A1} his definition into a systematic pathwise differentiation theory with respect to Brownian motion that is a natural counterpart to It\^o's Integral calculus; with a fundamental theorem of stochastic calculus relating this derivative to It\^o's integral, a differential stochastic chain rule, a differential stochastic mean value theorem, and other differentiation rules.  In \cite{A1,A2} it is shown that $S$ may be replaced with $f(S)$ for a reasonably large class of functions $f$.  We note briefly here that Allouba's definition of the stochastic quadratic covariation derivative (QCD) $\DW$ in \cite{A1}  is actually more general than \eqref{SAD}, enabling the differentiation in a more generalized sense, even when the derivative in \eqref{SAD} doesn't exist.  Namely,
\bdf[Allouba 2006 \cite{A1}: Definition 1.1]
\label{AlloubaDerivative} The stochastic difference and
stochastic derivative of a continuous semimartingale $S$ with respect to a Brownian motion $W$ are defined by%
\begin{equation}\lbl{ADisDer}
\mathcal{D}_{W_{t},h}S_{t}=\left\{ 
\begin{array}{cc}
\dsty\frac{3}{2h^{3}}\int_{0}^{h}r\left[ \left\langle S,W\right\rangle
_{t+r}-\left\langle S,W\right\rangle _{t-r}\right] dr; & 0<t<\infty ,\text{ }%
h>0 \\ 
\dsty\frac{3}{h^{3}}\int_{0}^{h}r\left\langle S,W\right\rangle _{r}dr; & t=0,%
\text{ }h>0%
\end{array}%
\right.
\end{equation}%
and
\begin{equation}\lbl{AD}
\mathbb{D}_{W_{t}}S_{t}=\lim_{h\rightarrow 0}\mathcal{D}_{W_{t},h}S_{t},
\end{equation}%
whenever this limit exists. If the derivative $\frac{d}{dt}\left\langle
S,W\right\rangle _{t}$ exists, then 
\begin{equation}\lbl{Astder}
\mathbb{D}_{W_{t}}S_{t}=\frac{d\left\langle S_{\cdot },W_{\cdot }\right\rangle _{t}}{dt},
\end{equation}
and $\DW$ is called the strong derivative of $S$ with respect to $W$.
The $k$-th $W$-derivative of $S$ is defined iteratively in the obvious way.  
\edf 
In \cite{Agn} and followup articles, the first author also generalizes his approach beyond the classical setting of Markov, semimartingale, or Gaussian processes to a much larger class of processes. 

For the rest of this paper, Let $W$ be a one-dimensional Brownian motion on the usual probability space $\OFFtP$ (the filtration satisfies the usual conditions of right continuity and completeness), where $\lbr\sFt\rbr$ is the augmentation under $\P$ of the natural filtration of $W$, $\lbr\sF^{W}_{t}\rbr_{t\in\Rp}$.  Let $T>0$ be arbitrary and fixed.  We denote by $\WrT$ the restriction of $W$ to the time interval $[0,T]$.   Unless stated otherwise,  our focus throughout this article will be on the strong derivative \eqref{Astder}.  Other distributional and Radon-Nikodym type versions of $\DW$---as well as obvious extensions to derivatives with respect to general semimartingales---and some of their implications are among many $\DW$-features discussed in \cite{A2}.  For more details on the quadratic covariation differentiation theory and its results, the reader is referred to Allouba's original article \cite{A1}.    

We first remark briefly on an aspect of that theory that is advantageous in our results here.  Since processes of bounded variations on compacts have quadratic variation zero, their QCD is identically $0$ (see Remark 1.1 in \cite{A1}, which says that these bounded variation processes are the ``constants'' in this quadratic covariation differential calculus).
An important consequential feature of the quadratic covariation derivative $\DW$ in \cite{A1} is that it is invariant under Girsanov's change of measure.  I.e., let $\Wt$ be the translated Brownian motion  $\Wt_{t}=W_{t}+\int_{0}^{t}\lambda \left( u\right) du\mbox{ for } \ 0\le t\le T$ and let $\Pt$ be the Girsanov changed probability measure, and assume the standard Novikov condition on $\lambda$ (see  Appendix \ref{appA} and \thmref{Girsanov} for the notation and setting and for a precise statement).  If $S$ is a continuous semimartingale and if either one of the two QCD derivatives $\DW S$ or $\mathbb{D}_{\widetilde{W}} S$ is finite  for $0\le t\le T$, then so is the other and they are indistinguishable ($\P$ and $\Pt$).  To see this, observe that
\begin{equation}\lbl{ADinv}
\begin{split}
\mathbb{D}_{\widetilde{W}_{t}}S_{t}&=\frac{d\left\langle S_{\cdot },\widetilde{W}_{\cdot }\right\rangle^{\Pt}_{t}}{dt}=\frac{d\left\langle S_{\cdot },W_{\cdot }+\int_{0}^{\cdot }\lambda \left(u\right) du\right\rangle^{\Pt} _{t}}{dt}\\
&=\frac{d\left\langle S_{\cdot },W_{\cdot }+\int_{0}^{\cdot }\lambda \left(u\right) du\right\rangle^{\P} _{t}}{dt}=\frac{d\left\langle S_{\cdot },W_{\cdot}\right\rangle^{\P} _{t}}{dt}=\mathbb{D}_{W_{t}}S_{t}; 
\end{split}
\end{equation}
for $\ t\in [0,T]$, a.s.~$\P$ and $\Pt$, where we used \lemref{invcomqcp} along with the fact that adding continuous processes of bounded variation on compacts does not alter the quadratic covariation process.  This invariance under change of measure feature results in a simpler representation in the $\DW$-variant of the Clark-Ocone formula under change of measure \eqref{ACOCOMconcl} than the classical one obtained using the Malliavin (or Hida-Malliavin) derivative \eqref{COCOMcncl} (see the original fundamental articles by Clark \cite{Cl}, Ocone \cite{Oc} and Ocone et al.~\cite{OcKaCOM,KOL}).  In fact, unlike the Malliavin or Hida Malliavin derivatives versions of the Clark-Ocone formula, the form of our variant in \thmref{ACO} under change of measure \eqref{ACOCOMconcl} is as simple as it is without change of measure \eqref{AFrep}.  This is true without requiring any further differentiability conditions on the Girsanov transform integrand beyond the standard Novikov condition for Girsanov theorem.   For a nice readable account and history of the Clark-Ocone formula in both the classical Malliavin and the Hida-Malliavin settings, we refer the reader to the excellent recent book by Di Nunno, \O ksendal, and Proske \cite{Levy} and the references therein.  For another  non-Malliavin and different Radon-Nikodym type approach we also refer the reader to Di Nunno's recent work \cite{DiN}.

The QCD $\mathbb{D}_{W}S=\left\{ \mathbb{D}_{W_{t}}S_{t};t\in [0,\infty)\right\}$, when it exists, is a stochastic process that is intimately connected to It\^{o}'s original construction of his stochastic integral via It\^{o}'s isometry using quadratic covariations; and it therefore leads to an approach to pathwise stochastic differentiation that is a natural  counterpart to It\^o's integration theory (see \cite{A1,A2} for more on this).   To wit, the derivative  $\mathbb{D}_{W}S$ is an anti-It\^{o}'s integral that yields a fundamental theorem of stochastic calculus (Theorem 2.1 and Theorem 2.2 in \cite{A1}), a differential stochastic mean value theorem (Lemma 2.1 in \cite{A1}), differential stochastic chain rules and more (Theorem 3.1 in \cite{A1} and also other versions in \cite{A2}).  In addition,  $\mathbb{D}_{W}S$ interacts with basic algebraic operations on semimartingales similarly to the action of the Newton elementary deterministic derivative on functions (Corollary 3.2 and Theorem 3.2 in \cite{A1}), making it a convenient tool for computations and proofs (\cite{A1,A2} and \thmref{ACO} below).  Several other extensions and applications, including a simple derivation of It\^{o}'s formula using this differentiation theory, are given in \cite{A1,A2}.  On the other hand, since $\mathbb{D}_{W}S$ in \eqref{SAD} (or \eqref{AD})  is a stochastic process defined in terms of quadratic covariations; it is a pathwise derivative that measures the rate of temporal change of a semimartingale $S$ (and reasonable functions thereof) with respect to temporal changes in a BM $W$ using the ``proper'' measure of time regularity of their H\"older-$\lpa1/2\rpa^{-}$ paths. This basic principle makes the differentiation theory in \cite{A1} amenable to generalizations that handle very general stochastic processes beyond the classical framework of Gaussian, Markovian, semimartingales processes.  This very general calculus theory is well beyond the scope of this article; and it is the subject of Allouba's program in \cite{Agn} and followup papers with Brownian-time processes (\cite{Abtp1,Abtp2}) and many other non classical processes.  

In this article we show that, even \emph{within} the It\^o setting, there are advantages to the stochastic differentiation theory in \cite{A1}.  Specifically, we apply it to derive and prove variants of the celebrated Clark-Ocone and Stroock formulas that are simple in form and proof (even under change of measure), and they are widely applicable. 
The proofs of our variants of the Clark-Ocone formula (\thmref{ACO} below) are simple consequences of the quadratic covariation differentiation theory in \cite{A1}---Theorems 2.1 and 2.2 in \cite{A1} (the QCD fundamental theorem of stochastic calculus), Theorem 3.1 and Corollary 3.1 in \cite{A1} (the QCD chain rules), and other QCD differentiation rules like Theorem 3.2 in \cite{A1}---along with It\^o's integral calculus. Since this stochastic differentiation theory in \cite{A1} is built using ingredients of It\^o's standard setup; the statement, proof, and applicability of \thmref{ACO} is naturally linked to It\^o's calculus setting, under $L^{2}$-type conditions, without the need for extra machinery and settings from distributional differentiation theory, Hida's white noise analysis, Malliavin or Hida-Malliavin calculi, and even without weakening the derivative to a Radon-Nikodym density.  The QCD Stroock variant (\thmref{AStr} below) is proved by an iterative application of our QCD Clark-Ocone formula. 
\subsection{The QCD variants of the Clark-Ocone formula with and without change of measure}
We denote by $F\in L^{2}\left(\Omega,\mathscr{F}_{T}, \mathbb{P}\right)$ an $L^{2}(\Omega,\P)$ and $\mathscr{F}_{T}$-measurable random variable.  The notion of almost indistinguishability in the sense of Theorem 2.2 \cite{A1} is useful for a more complete statement of our variants of the Clark-Ocone results.  
In the interest of moving quickly to the results, we refer the reader to Appendices \ref{appA} and \ref{B} for notations and for such definitions. 

The essence of our result is that, whether we change measure or not,  the integrand $X$ in the representation of an $L^{2}$-random variable $F$ is the stochastic process that is the derivative $\DW\E\lbk F|\sF\rbk$ (or $\mathbb{D}_{\widetilde{W}}\widetilde{\mathbb{E}}\left[ F|\mathscr{F}\right]$) of the naturally-associated martingale $\E\lbk F|\sF\rbk$\footnote{For any probability measure $\P$ defined on $\sFT\subset\sF$ and for any $Y\in L^{1}(\Omega,\P)$, we will always assume that $Y_{t}:=\E\lbk Y|\sFt\rbk$ is chosen from the equivalence class of $\E\lbk Y|\sFt\rbk$ in such a way that the resulting martingale $\E\lbk Y|\sFt\rbk=\lbr \E\lbk Y|\sFt\rbk;0\le t\le T\rbr$ has paths that are right continuous with left limits (RCLL or cadlag) almost surely.  This is of course possible by the right continuity and completeness of our filtration $\lbr\sFt\rbr$.  Of course, this also means that if $\E\lbk Y|\sFt\rbk$ is a modification of a continuous process $X$, then they are indistinguishable and $\E\lbk Y|\sFt\rbk$ is continuous almost surely.} (or $\Et\lbk F|\sF\rbk$) with respect to the BM $W$ (or $\Wt$), with respect to which we are integrating.  No $L^{p}$ conditions of any kind are assumed on the derivative $\DW$.
\begin{thm}[The QCD variants of the Clark-Ocone formula with and without change of measure]\label{ACO}$$ $$
\begin{enumerate}\renewcommand{\labelenumi}{$($\alph{enumi}$)$}
\item  
Assume that the random variable $F\in L^{2}\left( \Omega, \mathscr{F}_{T},\mathbb{P}\right)$.  
Then there exists an almost indistinguishable extension of $\DWEFsF$, $\DWEFsFaie$, such that $\DWEFsFaie\in\sPtwosrT$ and
\beq\lbl{Arepaie}
F=\mathbb{E}\left[ F\right] +\int_{0}^{T}\mathbb{D}^{\sc{aie}}_{W_{t}}\mathbb{E}\left[F|\mathscr{F}_{t}\right] dW_{t},\text{ a.s. }\mathbb{P}.
\eeq
If the process $\DWEFsF=\lbr\DWtEFsFt; t\in[0,T] \rbr$ is $\mathscr{B}([0,T])\times 
\mathscr{F}$ measurable; then $\DWEFsF\in\sPtworT$, and it is the unique---in the sense of almost indistinguishability \eqref{aiinteg}---process such that 
\begin{equation}\lbl{AFrep}
F=\mathbb{E}\left[ F\right] +\int_{0}^{T}\mathbb{D}_{W_{t}}\mathbb{E}\left[
F|\mathscr{F}_{t}\right] dW_{t},\text{ a.s. }\mathbb{P}.
\end{equation}
In particular, if $f:\R\to\R$ is either a bounded Borel-measurable function or a locally bounded Borel-measurable function with $\ds\lim_{x\to\pm\infty}x^{-2}\log^{+}\lab f(x)\rab^{2}=0,$ and if $F=f(W_{T});$ then the process $\DWEFsF$ is almost surely continuous and $F$ admits the representation \eqref{AFrep}.  
\item Suppose $F$ is $\mathscr{F}_{T}$ measurable.  Assume Novikov's condition $($\eqref{GirsanovCondition} $(ii)$ in \thmref{Girsanov}$)$ holds and assume that
\begin{equation}\lbl{ACOCOMcnd}
(i)\ \E\lbk Z_{T}^{2}F^{2}\rbk=\widetilde{\mathbb{E}}\left[Z_{T} F^{2}\right] <\infty \mbox{ and }(ii)\ \E F^{2}<\infty,
\end{equation}
where $Z_{T}$ is the Radon-Nikodym derivative in Girsanov's change of measure \thmref{Girsanov}.
Suppose further that the processes $\mathbb{D}_{\widetilde{W}}\widetilde{\mathbb{E}}\left[ F|\mathscr{F}\right]$ and $\DW\E\left[ Z_{T} F|\sF\right]$ are $\mathscr{B}([0,T])\times\mathscr{F}_{T}$ measurable, then 
\begin{equation}\lbl{ACOCOMconcl}
F=\widetilde{\mathbb{E}}\left[ F\right] +\int_{0}^{T}\mathbb{D}_{\widetilde{W%
}_{s}}\widetilde{\mathbb{E}}\left[ F|\mathscr{F}_{s}\right] d\widetilde{W}_{s},\text{ a.s. }\Pt\ (\mbox{and a.s. } \P) 
\end{equation}
If the measurability condition on $\DW\E\left[ Z_{T} F|\sF\right]$, and $\mathbb{D}_{\widetilde{W}}\widetilde{\mathbb{E}}\left[ F|\mathscr{F}\right]$ is dropped, then \eqref{ACOCOMconcl} holds with $\mathbb{D}_{\widetilde{W}}\widetilde{\mathbb{E}}\left[ F|\mathscr{F}\right]$ replaced by an almost indistinguishable extension $\mathbb{D}_{\widetilde{W}}^{\sc{aie}}\widetilde{\mathbb{E}}\left[ F|\mathscr{F}\right]$. 
\end{enumerate}
\end{thm}
Several observations are in order here and are summarized in the remarks below.
\brm
\begin{itemize}
\item We emphasize here that  the strong derivative $\DW$ from \cite{A1} that we use here in \thmref{ACO} is a derivative of  a function (for a.s.~$\omega$), defined in terms of the derivative of the quadratic covariation process with respect to time $t$, not a weak distributional type derivative or a Radon-Nikodym type derivative or its density (such weaker versions of $\DW$ and some of their implications are given in \cite{A2}).   This is an important feature of \thmref{ACO} since our variant of the Clark-Ocone formula is stated under conditions that are comparable to those that are given for the \emph{weak} Hida-Malliavin derivative (see Theorem 6.35 and Theorem 6.41 in \cite{Levy} which require the use of white noise analysis combined with Malliavin calculus as explained in Chapters 5 and 6 in \cite{Levy}).  These conditions allow us to handle many applications where the classical Malliavin differentiability condition $F\in\sDot$ (see \cite{Levy,Mal,Oc} and Appendix \ref{appB} below for Malliavin calculus background) is too strong of a condition, as we shall shortly see using an example from mathematical finance.   In this famous example, $F=\mathbb{I}_{[K,\infty)}\left( W_{T}\right)$ is the payoff of a digital option, where $K>0$ is a constant and $\mathbb{I}_{[K,\infty )}(\cdot)$ is the indicator function on the interval $[K,\infty)$.  It is well known that  $\mathbb{I}_{[K,\infty )}\left( W_{T}\right)\notin\sDot$  (see \cite{Levy}), but we show that it leads to a process $\E[\mathbb{I}_{[K,\infty )}\left( W_{T}\right)|\sF]$ that is \emph{infinitely} differentiable with respect to $\DW$ (see Subsection \ref{indicrepsec} below).   Also, we note that to obviate the need for the almost indistinguishable extension of $\DW$ we \emph{only} assume the measurability of $\DW$ (no $L^{2}$ conditions of any kind are assumed on $\DW$ since they follow for free as is clear from the proof below).  That measurability easily holds for a large class of random variables---including the Brownian indicator $\mathbb{I}_{[K,\infty)}\left( W_{T}\right)$---that are \emph{not} standard Malliavin differentiable (not in $\sDot$).                                                                                                               
\item The Clark-Ocone formula \cite{Cl,Oc} was extended by Karatzas and Ocone in \cite{OcKaCOM} to the Clark-Ocone formula under change of measure, a result that has proved very beneficial in mathematical finance.  Looking at the representation in the Clark-Ocone theorem under change of measure in the Malliavin setting  (e.g. \thmref{COCOM}),  \eqref{COCOMcncl}, we see that it is not as simple as its original version in \thmref{sCO} (this is true even when using the weak Hida-Malliavin derivative as in Theorem 6.35 and Theorem 6.41 in \cite{Levy}).  The QCD ($\DW$) variant, however, retains the simplicity of its representation, which has the same form in \thmref{ACO} (b) as it does in its unchanged measure version (\thmref{ACO} (a)).  It does so, without any differentiability requirements on the Girsanov transform integrand $\lambda$ beyond the conditions already demanded by Girsanov's \thmref{Girsanov} (this is not true even in the Hida-Malliavin setting see \cite{Levy} p.~107).   This simplicity is very useful in deriving the QCD Stroock's formula variant under change of measure, and many other examples, including when $\lambda =f\left( W\right)$ for $f$ that is only bounded and measurable.
\item   \thmref{ACO} (a) tells us that the integrand $X$ in It\^o's representation theorem is an almost indistinguishable version or extension of $\mathbb{D}_{W}\mathbb{E}\left[ F|\mathscr{F}\right]$, which is the quadratic covariation derivative  of the 
natural martingale associated with $F$$: \E\lbk F|\sF\rbk=\lbr\E\lbk F|\sFt\rbk,\sFt;0\le t\le T\rbr.$
Viewed this way, the integrand $X$ in It\^o's famous representation theorem is the stochastic process that is the rate of change of the martingale $\E\lbk F|\sF\rbk$ with respect to the BM with respect to which we are integrating.  
\item Note that the condition \eqref{ACOCOMcnd} together with H\"older inequality immediately imply \beq\lbl{tL2F}\Et F^{2}=\E Z_{T}F^{2}<\infty.\eeq
Another obvious but important remark here is that we cannot simply apply \thmref{ACO} (a) directly to the new Brownian motion $\Wt$ to get a representation with respect to $\Wt$, since $F$ is only assumed to be $\sFT$ measurable, and $\sFTt\subset\sFT$, where $\sFTt$ is the $\sigma$ algebra generated by $\lbr\Wt_{t};0\le t\le T\rbr$. 
\end{itemize}
\erm 
\subsection{The QCD Variant of Stroock's Formula with and without change of measure}
In his fundamental article \cite{Stroock}, Stroock identifies the integrands of the chaos expansion
of an $L^{2}$ random variable.  Using an iterated application of \thmref{ACO}, we get a $\DW$-variant of Stroock's formula.  We use the notations $J_{n}(g_{n})$ and $I_{n}(g_{n})$ for the $n$-fold iterated It\^o-Wiener integral over the simplex  $$\S_{n}=\left\{ \left( t_{1},t_{2},...,t_{n}\right);0\leq t_{1}\leq t_{2}\leq ...\leq t_{n}\leq T\right\}\subset[0,T]^{n}$$ and over $[0,T]^{n}$, respectively $($see \cite{Levy} and Appendix \ref{appB}$)$.   I.e.,
\beq\lbl{iwint}
\bsp
J_{n}\left( g_{n}\right) &=\int_{\S_{n}}g_{n}\left( t_{1},...,t_{n}\right) dW_{t_{1}}dW_{t_{2}}...dW_{t_{n-1}}dW_{t_{n}}\\
I_{n}\left( \hat{f}_{n}\right) &=\int_{\left[ 0,T\right] ^{n}}\hat{f}_{n}\left( t_{1},...,t_{n}\right) dW_{t_{1}}dW_{t_{2}}...dW_{t_{n-1}}dW_{t_{n}}
\end{split}
\eeq
for $g_{n}\in L^{2}\left(\S_{n}\right)$ and $\hat{f}_{n}\in \hat{L}^{2}\left([0,T] ^{n}\right)$, where $\hat{L}^{2}\left([0,T]^{n}\right) $ is the space of ${L}^{2}\left([0,T]^{n}\right)$ symmetric functions.
\begin{thm}[The QCD Variant of Stroock's Formula with and without change of measure]\label{AStr}$$ $$
\begin{enumerate}\renewcommand{\labelenumi}{$($\alph{enumi}$)$}
\item 
Suppose that $F\in L^{2}\OFTP$, with chaos expansion 
$$F=\sum_{n=0}^{\infty}J_{n}\lpa g_{n}\rpa=\sum_{n=0}^{\infty}I_{n}\left( \ft_{n}\right),$$  
where $J_{n}$ and $I_{n}$ are given by \eqref{iwint}. 
Assume that the random field $\vfi_{n}:\Omega\times S_{n}\to\R$ given by
\beqs
\bsp
\vfi_{n}(t_{1},t_{2},\ldots,t_{n})=\DWtone\E\lbk\DWtt\E\lbk\ldots\DWtn\E\lbk F\lgab\sFtn\rbk\cdots\lgab\sFtt\rbk\lgab\sFto\rbk;
 \end{split}
\eeqs
is almost surely continuous in $t_{1},t_{2},...,t_{n}$ for every $n=1, 2, \ldots$.
Let $$\Pi=\lbr\pi=(\pi_{1},\ldots,\pi_{n}); \pi \mbox{ is a permutation of }(1,\ldots,n),\ n\ge1\rbr.$$ Then, $J_{0}(g_{0})=g_{0}=\ft_{0}=I_{0}(\ft_{0})=\E F$; and for every $n\ge1$, every $0\leq t_{1}\leq t_{2}\leq ...\leq t_{n}\leq T$, and every one of the $n!$ permutations $(t_{\pi_{1}},\ldots,t_{\pi_{n}})$ the chaos expansion coefficients are given by
\begin{equation}\lbl{AStrcoef}
\bsp
\ft_{n}\left( t_{\pi_{1}},...,t_{\pi_{n}}\right)&=\frac{1}{n!}g_{n}(t_{1},\ldots,t_{n})
\\&=\frac{1}{n!}\E\lbk\DWtone\E\lbk\ldots\DWtn\E\lbk F\lgab\sFtn\rbk\cdots\lgab\sFto\rbk\rbk. 
\end{split}
\end{equation}
almost surely $\P$.
\item  Assume $F$ has the chaos expansion $F=\sum_{n=0}^{\infty }\widetilde{J}_{n}\left( g_{n}\right),$ where $\widetilde{J_{n}}$ denotes the $n$-fold iterated It\^o-Wiener integral with respect to $\widetilde{W}$ over the set $S_{n}$.  Assume the conditions of \thmref{ACO} (b) hold.  Suppose further that 
 the random fields $\vfi_{n}, \psi_{n}:\Omega\times S_{n}\to\R$ given by
\beq\lbl{Strcnd1}
\bsp
\vfi_{n}(t_{1},t_{2},\ldots,t_{n})&=\DWttone\Et\lbk\DWttt\Et\lbk\ldots\DWttn\Et\lbk F\lgab\sFtn\rbk\cdots\lgab\sFtt\rbk\lgab\sFto\rbk
\\\psi_{n}(t_{1},t_{2},\ldots,t_{n})&=\DWtone\E\lbk Z_{T}\DWttt\Et\lbk\ldots \DWttn \Et\lbk F\lgab\sFtn\rbk\cdots\lgab\sFtt\rbk\lgab\sFto\rbk
\end{split}
\eeq
are almost surely continuous in $t_{1},t_{2},...,t_{n}$ and that 
\beq\lbl{Strcnd2}
\bsp
\E Z_{T}^{2}\vfi_{n}^{2}(t_{1},t_{2},\ldots,t_{n})<\infty
\end{split}
\eeq for every $n=1, 2, \ldots$.Then, $g_{0}=\Et F$; and for every $n\ge1$ and every $0\leq t_{1}\leq t_{2}\leq ...\leq t_{n}\leq T$, the chaos expansion coefficients are 
\beq\lbl{AStrcoefCOM}
g_{n}(t_{1},\ldots,t_{n})=\Et\lbk\DWttone\Et\lbk\ldots\DWttn\Et\lbk F\lgab\sFtn\rbk\cdots\lgab\sFto\rbk\rbk\mbox{ a.s. }\P\mbox{ and }\Pt. 
\eeq
\end{enumerate}
\end{thm}
\brm\lbl{strundercomrem}
Part (b) in \thmref{AStr} makes it clear that another difference between here and the Malliavin setting is that here we have a simple Stroock's formula under change of measure by iterating \thmref{ACO} (b); whereas---because of the complexity of the representation \eqref{COCOMcncl} in the conclusion of the standard Clark-Ocone under change of measure (\thmref{COCOM})---this approach would quickly become significantly more complicated using either the Malliavin or Hida-Malliavin calculi.   
\erm
\section{Proofs of the QCD variants of Clark-Ocone and Stroock formulas}\label{pa}
\subsection{Proofs of the QCD variants of the Clark-Ocone formulas}\lbl{COsec}

We now give the proofs of our Clark-Ocone variants in \thmref{ACO}.  We start with the QCD variant under no change of measure.
\bpfs{Proof of \thmref{ACO} (a)}
Since $\E F^{2}<\infty$ and $\lbr\sFt\rbr$ is the augmented Brownian filtration (of $W$), then $\E\lbk F|\sF\rbk=\lbr\E\lbk F|\sFt\rbk,\sFt;0\le t\le T\rbr$ is a square-integrable RCLL $W$-Brownian martingale with $\E\lbk F|\sF_{0}\rbk=\E\lbk F\rbk$.  Therefore, by the Brownian martingale representation theorem (e.g., \cite{KaratzasShreve} p.~182)
\beq\lbl{Bmart}
\E\lbk F|\sFt\rbk=\E\lbk F\rbk+\int_{0}^{t}X_{s} dW_{s}; \quad 0\le t\le T,\mbox{ a.s. }\P,
\eeq
for some unique (in the sense of \eqref{aiinteg}) $X\in\sPtwosrT$.
Applying $\DW$ to both sides of \eqref{Bmart} and using Theorem 2.2 in \cite{A1} (the second QCD fundamental theorem of stochastic calculus)  yield a subset $\Omega^{*}\subset\Omega$, with $\P\lpa \Omega^{*}\rpa=1$; and a collection of zero Lebesgue-measure random sets $\lbr Z(\omega); \omega\in\Omega^{*}\rbr$ such that 
\beq\lbl{DWXconn}
\bsp
\DWtEFsFt\lpa\omega\rpa
=X_{t}\lpa\omega\rpa;\quad t\in[0,T]\setminus Z(\omega), \omega\in\Omega^{*}.
\end{split}
\eeq
I.e., $\DWEFsFaie:=X$ is an almost indistinguishable extension of $\DWEFsF$, and the representation in \eqref{Arepaie} is proved.  

Now, assume that $\DWEFsF:\Omega\times[0,T]\to\R$ is measurable.  To show the adaptability of $\mathbb{D}_{W}\mathbb{E}\left[ F|\mathscr{F}\right]$  ($(2)$ of  \defnref{intgrnds}), it is enough to show the adaptability of the right time derivative $\frac{{d}}{dt^{+}}\lqv\E\lbk F|\sF\rbk,W\rqv_{t}$.  But, $\frac{{d}}{dt^{+}}\lqv\E\lbk F|\sF\rbk,W\rqv_{t}\in\dsty\bigcap_{s>t}\sFs=\sFtp=\sFt$, where the last equality follows by the right continuity of the filtration $\lbr\sFt\rbr$.  This means the process $\mathbb{D}_{W}\mathbb{E}\left[ F|\mathscr{F}\right]$  satisfies the first two conditions in	\defnref{intgrnds}.      
This and \eqref{DWXconn} easily imply that \eqref{aiinteg} holds, with the $X$ in \eqref{Bmart} and $Y=\DWEFsF$, and so $X$ and $\DWEFsF$ are almost indistinguishable versions of each other and 
\begin{equation}\lbl{intXandDW}
\bsp
&\mathbb{E}\int_{0}^{T}\left\vert \mathbb{D}_{W_{t}}\mathbb{E}\left[ F|\mathscr{F}_{t}\right] \right\vert ^{2}dt<\infty
\lpa\mbox{hence $\DW\E\lbk F|\sF\rbk\in\sPtworT$}\rpa, \mbox{ and}\\
&\E\lbk F|\sFt\rbk=\E\lbk F\rbk+\int_{0}^{t}\DWs\E\lbk F|\sFs\rbk dW_{s}; \quad 0\le t\le T\mbox{ a.s. }\P.
\end{split}
\end{equation}
We are done by setting $t=T$ in \eqref{intXandDW} since $\E\lbk F|\sFT\rbk=F$.  

Finally, if $f:\R\to\R$ is either a bounded Borel-measurable function or a locally bounded Borel-measurable function with $\ds\lim_{x\to\pm\infty}x^{-2}\log^{+}\lab f(x)\rab^{2}=0$ then $F=f(W_{T})\in L^{2}(\Omega,\sFT,\P)$ and $v(T-t,W_{t})=\E[f\left( W_{T}\right)|\sFt]$ is $\mathrm{C}^{1,2}$ (e.g., \cite{Dur} pp.~128--130); and so It\^o's rule followed by Theorem 2.1 in \cite{A1} (the QCD fundamental theorem of stochastic calculus) implies that $\mathbb{D}_{W_{t}} \E[f\left( W_{T}\right)|\sFt]=\partial_{2}v(T-t,W_{t})$ (the first partial derivative in the second variable evaluated at $W_{t}$) is continuous in $t$ almost surely, and $F=f(W_{T})$ admits the representation \eqref{AFrep}. 
\epfs
Next, we prove the QCD variant of Clark-Ocone under change of measure.  
\bpfs{Proof of \thmref{ACO} $(b)$}
Let 
\begin{equation}\lbl{Y00}
Y_{t}:=\widetilde{\mathbb{E}}\left[ F|\mathscr{F}_{t}\right]
\end{equation}%
and notice that  \eqref{tL2F} and Jensen's inequality yield 
\beq\lbl{tL2Yt
}\Et Y_{t}^{2}=\Et\lpa \widetilde{\mathbb{E}}\left[ F|\mathscr{F}_{t}\right]\rpa^{2}\le \Et\lpa \widetilde{\mathbb{E}}\left[ F^{2}|\mathscr{F}_{t}\right]\rpa=\Et F^{2}<\infty.
\eeq
Let
\begin{equation}\lbl{ZLm}
\begin{split}
\Lambda_{t} 
=Z^{-1}_{t}=\exp \left[ \int_{0}^{t}\lambda _{s}d\widetilde{W}%
_{s}-\frac{1}{2}\int_{0}^{t}\lambda _{s}^{2}ds\right],
\end{split}
\end{equation}%
where $Z$, $\lambda$, and $\widetilde{W}$ are as in the change of measure \thmref{Girsanov}.  Now applying 
\lemref{CondExpCOM} and \thmref{ACO} (a) to $\mathbb{E}\left[ Z_{T} F|\mathscr{F}_{t}\right]$, using \eqref{ACOCOMcnd}, we obtain
\begin{equation}\lbl{Y0}
\begin{split}
Y_{t}&=\Lambda_{t}  \mathbb{E}\left[ Z_{T} F|\mathscr{F}_{t}\right] 
=\Lambda_{t} \left[ \mathbb{E}\left[ Z_{T} F\right]+\int_{0}^{t}\mathbb{D}_{W_{s}}\mathbb{E}\left[ Z_{T}F|%
\mathscr{F}_{s}\right] dW_{s}\right]
\\&=:\Lambda_{t} U_{t}; \ 0\le t\le T\mbox{ a.s. }\P\mbox{ and }\Pt.
\end{split}
\end{equation}
It\^{o}'s formula easily gives $d\Lambda_{t} =\Lambda_{t} \lambda_{t}d%
\widetilde{W}_{t}$, and integration by parts for It\^o's calculus then gives
\begin{equation}\lbl{Y1}
\begin{split}
dY_{t}&=\Lambda_{t} dU_{t} +U_{t} d\Lambda_{t} +d\left\langle \Lambda ,U\right\rangle _{t}
\\&=\Lambda_{t} \mathbb{D}_{W_{t}}\mathbb{E}\left[Z_{T}
F|\mathscr{F}_{t}\right] dW_{t}+U_{t} \Lambda_{t} \lambda \left(
t\right) d\widetilde{W}_{t}
+\Lambda_{t} \lambda_{t}\mathbb{D}_{W_{t}}\mathbb{E}%
\left[ Z_{T} F|\mathscr{F}_{t}\right] dt
\\ &=\lbk\Lambda_{t} \mathbb{D}_{W_{t}}\mathbb{E}\left[ Z_{T}F|\mathscr{F}_{t}\right] +\lambda_{t}Y_{t}\rbk d\widetilde{W}_{t}.
\end{split}
\end{equation}
  On the other hand, by Theorem 2.2 in \cite{A1} (the QCD fundamental theorem of stochastic calculus), Theorem 3.2 in \cite{A1} (the QCD product rule), and the QCD invariance under change of measure \eqref{ADinv}, we obtain  that $\exists$ a set $\Omega^{*}\subset\Omega$ $\ni$ $\mathbb{P}(\Omega^{*})=1$, and for each $\omega\in\Omega^{*}$ there  is a Lebesgue-measure-zero random set $Z(\omega)\subset[0,T]$ such that
\begin{equation}\lbl{Y2}
\begin{split}
&\mathbb{D}_{W_{t}}\left( \Lambda_{t} \mathbb{E}\left[ Z_{T} F|\mathscr{F}_{t}\right] \right) 
=\lambda_{t}\Lambda_{t}\mathbb{E}\left[ Z_{T} F|\mathscr{F}_{t}\right] +\Lambda_{t} \mathbb{D}_{W_{t}}\mathbb{E}\left[
Z_{T} F|\mathscr{F}_{t}\right] 
\\&=\lambda_{t}Y_{t} +\Lambda_{t} \mathbb{D}
_{W_{t}}\mathbb{E}\left[ Z_{T} F|\mathscr{F}_{t}\right];\   t\in \left[0,T\right]\setminus Z(\omega),\ \omega\in\Omega^{*}, 
\end{split}
\end{equation}
where we used a trivially obvious adaptation of the proof of Theorem 3.2 in \cite{A1} to account for the possibly-discontinuous integrand case.  But equations \eqref{Y00} and \eqref{Y0} in conjunction with the $\DW$-invariance under change of measure \eqref{ADinv} imply the indistinguishability 
\beq\lbl{the2deriv}
\bsp
\mathbb{D}_{\widetilde{W}}\widetilde{\mathbb{E}}\left[ F|\mathscr{F}\right]=\DWtld \Lambda\E\left[ Z_{T}F|\sF\right]=\DW \Lambda\E\left[ Z_{T}F|\sF\right]; \mbox{ a.s.~}\P\mbox{ and }\Pt.
\end{split}
\eeq
The measurability of $\DW  \Lambda\E\left[ Z_{T}F|\sF\right]$ follows from \eqref{the2deriv} together with the measurability assumption on $\mathbb{D}_{\widetilde{W}}\widetilde{\mathbb{E}}\left[ F|\mathscr{F}\right]$.   Using \eqref{Y1}, \eqref{Y2}, and \eqref{the2deriv} we then have
\begin{equation}\lbl{Y3}
\begin{split}
dY_{t}&=\mathbb{D}_{W_{t}}\left( \Lambda_{t} \mathbb{E}\left[ Z_{T} F|\mathscr{F}_{t}\right] \right) d\widetilde{W}_{t}=\mathbb{D}_{\widetilde{W}_{t}}
\widetilde{\mathbb{E}}\left[ F|\mathscr{F}_{t}\right] d\widetilde{W}_{t};  \mbox{ a.s.~}\P\mbox{ and }\Pt.
\end{split}
\end{equation}%
By \eqref{Y00} we see that $Y_{T} =F$ and $Y_{0} =\widetilde{\mathbb{E}}%
\left[ F\right]$, and the desired conclusion follows.

Finally, if the measurability assumption on $\mathbb{D}_{\widetilde{W}}\widetilde{\mathbb{E}}\left[ F|\mathscr{F}\right]$ and $\DW\E\left[ Z_{T} F|\sF\right]$ is dropped;
then by \thmref{ACO} (a) and \eqref{ACOCOMcnd}---together with the argument after \eqref{Y0}---equations \eqref{Y0} through \eqref{Y2} hold with an almost indistinguishable extension $\DW^{\sc{aie}}\E\left[ Z_{T}F|\sF\right]\in\sPtwosrT$ in place of $\DW\E\left[ Z_{T} F|\sF\right]$.  In particular,
\begin{equation}\lbl{Y2aie}
\begin{split}
&\mathbb{D}_{W_{t}}\left( \Lambda_{t} \mathbb{E}\left[ Z_{T}F|\mathscr{F}_{t}\right] \right) 
=\lambda_{t}Y_{t} +\Lambda_{t} \DWt^{\sc{aie}}\E\left[ Z_{T}F|\sFt\right];\ t\in \left[0,T\right]\setminus Z(\omega),\ \omega\in\Omega^{*}, 
\end{split}
\end{equation}
  Now, the process on the  right hand side of \eqref{Y2aie}
\beq\lbl{rhsprocess}
\lbr\lambda_{t}(\omega)Y_{t}(\omega) +\Lambda_{t}(\omega) \DWt^{\sc{aie}}\E\left[ Z_{T}F|\sFt\right](\omega); (t,\omega)\in[0,T]\times\Omega\rbr
\eeq  
defines an almost indistinguishable extension of  
$\mathbb{D}_{W}\left( \Lambda\mathbb{E}\left[ Z_{T}F|\mathscr{F}_{t}\right] \right)$, which we call $\DW^{\sc{aie}} Z^{-1}\E\left[ Z_{T} F|\sF\right]$. On the other hand, the indistinguishability of $\widetilde{\mathbb{E}}\left[ F|\mathscr{F}\right]$ and $\Lambda \mathbb{E}\left[ Z_{T} F|\mathscr{F}\right] $ in equation \eqref{Y0} together with \eqref{Y2aie} and the $\DW$-invariance under change of measure \eqref{ADinv} imply that 
\begin{equation}\lbl{Y17aie}
\begin{split}
\mathbb{D}_{\widetilde{W}_{t}}\widetilde{\mathbb{E}}\left[ F|\mathscr{F}_{t}\right]
=\lambda_{t}Y_{t} +\Lambda_{t} \DWt^{\sc{aie}}\E\left[ Z_{T}F|\sFt\right];\ t\in \left[0,T\right]\setminus Z(\omega),\ \omega\in\tilde{\Omega}, 
\end{split}
\end{equation}
for some $\tilde{\Omega}\subset\Omega$ with $\P\lpa\tilde{\Omega}\rpa=\Pt\lpa\tilde{\Omega}\rpa=1$.  This  means that the  process in \eqref{rhsprocess} is also an almost indistinguishable extension of $\mathbb{D}_{\widetilde{W}}\widetilde{\mathbb{E}}\left[ F|\mathscr{F}\right]$, call it $\mathbb{D}_{\widetilde{W}}^{\sc{aie}}\widetilde{\mathbb{E}}\left[ F|\mathscr{F}\right]   $; in particular, we have the indistinguishability
$$\DW^{\sc{aie}} Z^{-1}\E\left[ Z_{T} F|\sF\right]=\lambda Y +\Lambda \DW^{\sc{aie}}\E\left[ Z_{T}F|\sF\right]=\mathbb{D}_{\widetilde{W}}^{\sc{aie}}\widetilde{\mathbb{E}}\left[ F|\mathscr{F}\right];  \mbox{ a.s.~}\P\mbox{ and }\Pt,$$
and we have 
\begin{equation}\lbl{Y3aie}
\begin{split}
dY_{t}&=\DWt^{\sc{aie}} Z^{-1}_{t}\E\left[ Z_{T} F|\sF\right] d\widetilde{W}_{t}=\mathbb{D}_{\widetilde{W}_{t}}^{\sc{aie}}\widetilde{\mathbb{E}}\left[ F|\mathscr{F}\right]d\widetilde{W}_{t};  \mbox{ a.s.~}\P\mbox{ and }\Pt.
\end{split}
\end{equation} 
The proof is complete. 
\epfs
\subsection{Proofs of the QCD variants of Stroock's formulas}
We now turn to the proof of the QCD variant of Stroock's formula.  Since the proof reduces to a simple iteration of \thmref{ACO} along with an adaptation of the standard chaos expansion proof (e.g., see \cite{Levy}), we simply indicate the changes.
We now give the proof of the version under no change of measure. 

\bpfs{Proof of \thmref{AStr} $(a)$}
By \thmref{ACO} (a),  
\begin{equation}\lbl{Str1}
F=\mathbb{E}\left[ F\right] +\int_{0}^{T}\mathbb{D}_{W_{s_{1}}}\mathbb{E}\left[ F|\mathscr{F}_{s_{1}}\right] dW_{s_{1}}.
\end{equation}%
Applying \thmref{ACO} (a) to $\mathbb{D}_{W_{s_{1}}}\mathbb{E}\left[ F|\mathscr{F}_{s_{1}}\right]$, $0\le s_{1}\le T$, we have that 
\begin{equation}\lbl{Str2}
\bsp
&\mathbb{D}_{W_{s_{1}}}\mathbb{E}\left[ F|\mathscr{F}_{s_{1}}\right] 
\\&=\mathbb{E}\left[ \mathbb{D}_{W_{s_{1}}}\mathbb{E}\left[ F|\mathscr{F}_{s_{1}}\right] \right] +\int_{0}^{s_{1}}\mathbb{D}_{W_{s_{2}}}\mathbb{E}\left[ \mathbb{D}_{W_{s_{1}}}\mathbb{E}\left[ F|\mathscr{F}_{s_{1}}\right]\lgab\mathscr{F}_{s_{2}}\right] dW_{s_{2}}.
\end{split}
\end{equation}%
Now, if we define $g_{0}=\E\lbk F\rbk$, $\vfi_{1}(s_{1})=\mathbb{D}_{W_{s_{1}}}\mathbb{E}\left[ F|\mathscr{F}_{s_{1}}\right]$, $g_{1}(s_{1})=\E\vfi_{1}(s_{1})$, and $\vfi_{2}(s_{2},s_{1})=\mathbb{D}_{W_{s_{2}}}\mathbb{E}\left[ \mathbb{D}_{W_{s_{1}}}\mathbb{E}\left[ F|\mathscr{F}_{s_{1}}\right]\lgab\mathscr{F}_{s_{2}}\right]$, for $0\le s_{2}\le s_{1}\le T$, and then substitute \eqref{Str2} into \eqref{Str1} we get
\beq\lbl{g0g1}
\bsp
F&=\mathbb{E}\left[ F\right] +\int_{0}^{T}\mathbb{E}\left[ \mathbb{D}_{W_{s_{1}}}\mathbb{E}\left[ F|\mathscr{F}_{s_{1}}\right] \right] dW_{s_{1}} 
\\&+\int_{0}^{T}\int_{0}^{s_{1}}\mathbb{D}_{W_{s_{2}}}\mathbb{E}\left[ \mathbb{D}_{W_{s_{1}}}\mathbb{E}\left[ F|\mathscr{F}_{s_{1}}\right]\lgab\mathscr{F}_{s_{2}}\right] dW_{s_{2}}dW_{s_{1}}
\\&=g_{0}+\int_{0}^{T}g_{1}(s_{1})dW_{s_{1}}+\int_{0}^{T}\int_{0}^{s_{1}}\vfi_{2}(s_{2},s_{1}) dW_{s_{2}}dW_{s_{1}} 
\end{split}
\eeq
Iterating this procedure we obtain, after $n$ steps, 
\beq\lbl{nit}
\bsp
F=\sum_{k=0}^{n}J_{k}(g_{k})+\int_{S_{n+1}}\vfi_{n+1}dW^{\otimes(n+1)},
\end{split}
\eeq
where $S_{n+1}=\left\{ \left( t_{1},t_{2},...,t_{n+1}\right);0\leq t_{1}\leq t_{2}\leq ...\leq t_{n+1}\leq T\right\}$, $n=0, 1, 2, \ldots$ and 
\beq\lbl{nitexp} 
\bsp
&\vfi_{n+1}(t_{1},t_{2},\ldots,t_{n+1})=\DWtone\E\lbk\DWtt\E\lbk\ldots\DWtnpo\E\lbk F\lgab\sFtnpo\rbk\cdots\lgab\sFtt\rbk\lgab\sFto\rbk,
\\&\int_{S_{n+1}}\vfi_{n+1}dW^{\otimes(n+1)}=\int_{0}^{T}\int_{0}^{\tnpo}\cdots\int_{0}^{\tt}\vfi_{n+1}(t_{1},t_{2},\ldots,t_{n+1})dW_{\tone}\cdots dW_{\tnpo},
\\&J_{k}(g_{k})=\int_{S_{k}}g_{k}dW^{\otimes(k)}; \ J_{0}(g_{0})=g_{0}=\E F,\ g_{k}=\E\vfi_{k}:S_{k}\to\R,\ k=1, 2,\ldots, n.
\end{split}
\eeq 

Now, letting $\Psi_{n+1}=\int_{S_{n+1}}\vfi_{n+1}dW^{\otimes(n+1)}$ for $n=0,1,2,\ldots$; it is clear that 
\beq\lbl{L2bdd}
\lnrm \Psi_{n+1}\rnrm^{2}_{L^{2}(\Omega,\P)}=\E\lbk\Psi^{2}_{n+1}\rbk\le \E F^{2}=\lnrm F\rnrm^{2}_{L^{2}(\Omega,\P)}<\infty; \ \forall n=0,1,2,\ldots
\eeq 
and It\^o's isometry implies that
\beq\lbl{L2Fassum}
\lnrm F\rnrm^{2}_{L^{2}(\Omega,\P)}=\sum_{k=0}^{n}\lnrm J_{k}(g_{k})\rnrm^{2}_{L^{2}(\Omega,\P)}+\lnrm \Psi_{n+1}\rnrm^{2}_{L^{2}(\Omega,\P)}
\eeq
and therefore $\sum_{k=0}^{\infty}J_{k}(g_{k})$ converges in $L^{2}(\Omega,\P)$; i.e., $\lim_{n\to\infty}\Psi_{n+1}=\Psi<\infty$ in $L^{2}(\Omega,\P)$.  Now, using the identical argument as in the standard chaos expansion proof (e.g., p.~14 in \cite{Levy}) it follows that $\Psi=0$.  Hence, we conclude that 
\begin{equation}\lbl{repovrS}
\bsp
&F=\sum_{k=0}^{\infty }J_{k}\lpa g_{k}\rpa=g_{0}+\sum_{k=1}^{\infty }\int_{S_{k}}g_{k}dW^{\otimes(k)}
=\E\lbk F\rbk\\&+\sum_{k=1}^{\infty }\int_{0}^{T}\int_{0}^{\tk}\cdots\int_{0}^{\tt}\E\lbk\DWtone\E\lbk\ldots\DWtk\E\lbk F\lgab\sFtk\rbk\cdots\lgab\sFto\rbk\rbk dW_{\tone}\cdots dW_{\tk},
\\&\lnrm F\rnrm^{2}_{L^{2}(\Omega,\P)}=\sum_{k=0}^{\infty}\lnrm J_{k}(g_{k})\rnrm^{2}_{L^{2}(\Omega,\P)}
\end{split}
\end{equation}
The only thing left is to rewrite the $n$-fold It\^o integral $J_{k}\lpa g_{k}\rpa=\int_{S_{k}}g_{k}dW^{\otimes(k)}$ over the subset $S_{k}\subset[0,T]^{k}$ as an $n$-fold It\^o integral $I_{k}\lpa \ft_{k}\rpa=\int_{[0,T]^{k}}\ft_{k}dW^{\otimes(k)}$ of appropriately related functions $\ft_{k}$ over $[0,T]^{k}$.  We start by extending $g_{k}$'s domain from $S_{k}$ to $[0,T]^{k}$ by setting 
$$g_{k}(t_{1},\ldots,t_{k})=0,\ (t_{1},\ldots,t_{n})\in[0,T]^{n}\setminus S_{n}$$
We then define $\ft_{k}$ to be the symmetrization of $g_{k}$:
$$\ft_{k}(t_{1},\ldots,t_{k}):=\frac{1}{k!}\sum_{\pi}g_{k}(t_{\pi_{1}},\ldots,t_{\pi_{k}});\ (\tone,\ldots,\tk)\in[0,T]^{k},$$
where the sum is taken over all $k!$ permutations $\pi$ of $(1,\ldots,k)$.  Then
\begin{equation*}
I_{k}\lpa\ft_{k}\rpa=k!J_{k}\left(\ft_{k}\right) =J_{k}\left( g_{k}\right), 
\end{equation*}%
completing the proof.
\epfs
\bpfs{Proof of \thmref{AStr} $(b)$}  The proof in this case follow the same iterative scheme of part (a) with obvious notational changes and using the assumed conditions along with  \thmref{ACO} (b) in place of \thmref{ACO} (a).  We omit the details.
\epfs

\section{Three applications for an $F$\ not in Malliavin's $\sDot$ space}
\subsection{Brownian Indicator representation}\lbl{indicrepsec}
An easy application of \thmref{ACO} above together with Theorem 3.1 in \cite{A1} (the QCD chain rule)  give us the following \emph{very} short proof of the representation of the Brownian indicator $F=\mathbb{I}_{[K,\infty )}\left( W_{T}\right)\notin\sDot$, for which the standard Malliavin derivative Clark-Ocone theorem doesn't apply, without the need for the many added technical aspects (including the use of Donsker's delta function) in the Hida-Malliavin derivative setting (see \cite{Levy} for a Hida-Malliavin derivation).
\bpr\label{PropositionIndicator}
If $K\in\R$ is fixed but arbitrary and if $p(t,x)=\frac{1}{\sqrt{2\pi t}}e^{-\frac{x^{2}}{2t}}$, then
\begin{equation}\lbl{indicrep}
\bsp
\mathbb{I}_{[K,\infty )}\left( W_{T}\right)
=\mathbb{P}\left[ W_{T}>K\right]+\int_{0}^{T}p(T-s,W_{s}-K)dW_{s}; \mbox{ a.s. }\P.
\end{split}
\end{equation}
\epr
\brm\lbl{appACO}
This application of \thmref{ACO} is an example of the fact that, in many cases, the measurability assumption on $\DWEFsF$ is not an added restriction.  Even for many $F$ that are not (standard) Malliavin differentiable, $\DWEFsF$ satisfies much more than this measurability condition.
\erm

\bpf
Let $F=\mathbb{I}_{[K,\infty )}\left( W_{T}\right)$.  By the Markov property and a simple change of variable we have for $0<s<T$
\begin{equation}\lbl{BIcndexpect}
\begin{split}
\mathbb{E}\left[ \mathbb{I}_{[K,\infty )}\left(W_{T}\right)|\mathscr{F}_{s}\right]&=\int_{-\infty}^{\infty}p(T-s,W_{s},y)\mathbb{I}_{[K,\infty )}\left(y\right)dy
\\&=\frac{1}{\sqrt{2\pi}}\int_{-\infty}^{\frac{W_{s}-K}{\sqrt{T-s}}}e^{-\frac{u^{2}}{2}}du
\end{split}
\end{equation}
By Theorem 3.1 in \cite{A1} (the QCD stochastic chain rule), 
\begin{equation}\lbl{DWofindic}
\bsp
\mathbb{D}_{W_{s}} \mathbb{E}\left[ \mathbb{I}_{[K,\infty )}\left(W_{T}\right)|\mathscr{F}_{s}\right]&=\frac{1}{\sqrt{2\pi \left(
T-s\right) }}\exp \left( \frac{-\left( W_{s}-K\right) ^{2}}{2\left(T-s\right) }\right)\\&=p(T-s,W_{s}-K);\  0<s<T \mbox{ a.s.}
\end{split}
\end{equation}%
Finally, $\mathbb{E}\left[ \mathbb{I}_{[K,\infty )}\left( W_{T}\right) \right] =\mathbb{P}\left[ W_{T}>K\right]$,
and the conclusion follows by \thmref{ACO} (a).
\epf
\brm\lbl{infdiffDW}
It easily follows from the quadratic covariation differentiation theory in \cite{A1} (Theorem 2.1 and Theorem 3.1 in \cite{A1}) that the conditional expectation $\mathbb{E}\left[ \mathbb{I}_{[K,\infty )}\left(W_{T}\right)|\mathscr{F}\right]$ is infinitely $\DW$-differentiable and that the $n$-th quadratic covariation derivative is given by
\beq\lbl{DWnthder}
\DWt^{(n)}\mathbb{E}\left[ \mathbb{I}_{[K,\infty )}\left(W_{T}\right)|\mathscr{F}_{t}\right]=
\partial_{2}^{n-1}p\lpa T-t,W_{t}-K \rpa;\  0<t<T, \ \mbox{a.s.}
\eeq
for $ n =1,2,\ldots$
\erm
\subsection{Application to digital options} \label{BinaryOptionSection}  As an example financial mathematics application of \thmref{ACO} (b) together with Theorem 3.1 in \cite{A1} (the QCD chain rule), we will find the replicating portfolio for a digital (binary) option  in the Black-Scholes-Merton framework. I.e., we assume a payoff of the form 
\begin{equation}\lbl{digitalpayoff}
V_{T}=\mathbb{I}_{[K,\infty )}\left( W_{T}\right) ,\text{ where }K>0\text{
is fixed,}
\end{equation}%
and where $V_{T}$ is the payoff at a fixed but arbitrary time of maturity $T$.  To carry out this analysis in the Malliavin derivative framework, some modifications are needed first;  e.g., working with the Hida-Malliavin derivative in the white noise setting (see \cite{Levy}).  Following the derivation in Shreve's book \cite{ShreveFinance}, assume that we have a stock whose price $P$ satisfies the SDE%
\begin{equation*}
dP_{t}=b_{t}P_{t}dt+a_{t}P_{t}dW_{t}.
\end{equation*}
where $b$ and $a$ are assumed deterministic and continuous.
Let $V_{t}$ and $X_{t}$ be $\mathscr{F}_{t}-$measurable random variables where $V_{t}$ is the payoff at time $t$ of
a derivative security and $X_{t}$ the portfolio value at time $t$. Our goal
is to find initial capital $X_{0}$ and number of shares invested at time $t,$
the portfolio $\Delta _{t},$ so that%
\begin{equation}\lbl{hedge}
D_{T}X_{T}=D_{T}V_{T}\text{ a.s.},
\end{equation}%
where the discount factor $D$ is defined as $D_{t}=e^{-\int_{0}^{t}r(s)ds}$, $r$ is the rate of return which we assume deterministic for simplicity,  and $V_{T}$ is the digital payoff in \eqref{digitalpayoff}.
The value of the portfolio at time $t$ is (see page 154 of \cite%
{ShreveFinance})%
\begin{equation*}
dX_{t}=\left[ r_{t}X_{t}+\left( b_{t}-r_{t}\right) \Delta _{t}P_{t}\right]
dt+a_{t}\Delta _{t}P_{t}dW_{t}.
\end{equation*}%
Assume that the market price of risk 
$
\lambda _{t}=\frac{b_{t}-r_{t}}{a_{t}}
$
satisfies Novikov's integrability condition (\eqref{GirsanovCondition} (ii)).  Then, under the usual risk-neutral measure $\widetilde{\mathbb{P}}$ in \thmref{Girsanov}, the process $DX$ given by 
\begin{equation}\label{CompletenessEquation}
D_{t}X_{t}=D_{0}X_{0}+\int_{0}^{t}\Delta _{s}a_{s}D_{s}P_{s}d\widetilde{W}_{s}; \ 0\le t\le T  
\end{equation}
is a martingale; i.e. 
\begin{equation*}
\widetilde{\mathbb{E}}\left[ D_{T}V_{T}|\mathscr{F}_{t}\right] =\widetilde{\mathbb{E}}\left[ D_{T}X_{T}|%
\mathscr{F}%
_{t}\right] =D_{t}X_{t}\text{ \ \ }t\in \lbrack 0,T],
\end{equation*}%
where $\widetilde{\mathbb{E}}$ is the expectation under $\widetilde{\mathbb{P}}$. Thus, as in Shreve \cite{ShreveFinance}, with $V_{t}$ being the price of the derivative security at time $t$ we get the risk neutral pricing formula
\begin{equation*}
D_{t}V_{t}=\widetilde{\mathbb{E}}\left[ D_{T}V_{T}|%
\mathscr{F}%
_{t}\right]; \ 0\le t\le T.
\end{equation*}%
Clearly, this means that $\lbr D_{t} V_{t},\sFt;0\le t\le T\rbr$ is a martingale under $\widetilde{\mathbb{P}}$.   By \thmref{ACO} (b) and by the deterministic assumption on $D$ and the fact that $D_{T}$ is independent of $t$, we get 
\begin{equation*}
D_{T}V_{T}=\widetilde{\mathbb{E}}\left[ D_{T}V_{T}\right] +\int_{0}^{T}D_{T}%
\mathbb{D}_{\widetilde{W}_{t}}\widetilde{\mathbb{E}}\left[ V_{T}|%
\mathscr{F}%
_{t}\right] d\widetilde{W}_{t}.
\end{equation*}%
Equation \eqref{hedge} will be satisfied if
\begin{equation*}
D_{T}\mathbb{D}_{\widetilde{W}_{t}}\widetilde{\mathbb{E}}\left[ V_{T}|\mathscr{F}_{t}\right] =\Delta _{t}a_{t}D_{t}P_{t},\ 0\le t\le T\mbox{ a.s.}
\end{equation*}
We now proceed to express the number of shares $\Delta_{t}$ at any time $t$ in terms of $a$, $\lambda$ and $W$.  To further simplify the computation  we assume $a\neq0$.  Now, the above discussion leads to 
\begin{equation*}
\Delta _{t}=e^{-\int_{t}^{T}r(s)ds}a_{t}^{-1}P_{t}^{-1}\mathbb{D}_{%
\widetilde{W}_{t}}\widetilde{\mathbb{E}}\left[ \mathbb{I}_{[K,\infty
)}\left( W_{T}\right) |%
\mathscr{F}%
_{t}\right] .
\end{equation*}

We could figure out the conditional expectation using the short steps in the proof of \propref{PropositionIndicator}; instead, we give a more financial mathematics argument by using the same idea used to derive the
solution of the Black-Scholes-Merton model in \cite{ShreveFinance}. Remembering that $\lambda$ is deterministic, we use
the Markov property to write 
\begin{equation*}
\begin{split}%
\widetilde{c}\left( t,W_{t}\right) 
&%
=\widetilde{\mathbb{E}}\left[ \mathbb{I}_{[K,\infty )}\left( W_{T}\right) |%
\mathscr{F}%
_{t}\right] =\widetilde{\mathbb{E}}\left[ \left.\mathbb{I}_{[K,\infty )}\left( 
\widetilde{W}_{T}-\int_{0}^{T}\lambda _{s}ds\right) \right|\mathscr{F}_{t}\right] 
\\ %
&%
=\widetilde{\mathbb{E}}\left[ \left.\mathbb{I}_{[K,\infty )}\left( \frac{%
\widetilde{W}_{T}-\widetilde{W}_{t}}{\sqrt{T-t}}\sqrt{T-t}+{W}_{t}-\int_{t}^{T}\lambda _{s}ds\right) \right|\mathscr{F}_{t}\right] 
\\ %
\end{split}%
\end{equation*}%
Now, $Y:=\frac{\widetilde{W}_{T}-\widetilde{W}_{t}}{\sqrt{T-t}}$ is a standard Normal random variable independent of $\mathscr{F}_{t},$ so that
\begin{equation*}
\widetilde{c}\left( t,x\right) =\int_{-\infty }^{\infty }\frac{1}{\sqrt{2\pi 
}}e^{-\frac{1}{2}y^{2}}\mathbb{I}_{[K,\infty )}\left( -y\sqrt{T-t}%
+x-\int_{0}^{T}\lambda _{s}ds\right) dy.
\end{equation*}%
Notice that $-y\sqrt{T-t}+x-\int_{t}^{T}\lambda _{s}ds>K$ if $y<\frac{%
x-\int_{t}^{T}\lambda _{s}ds-K}{\sqrt{T-t}}:=\widetilde{d}_{+}\left( x,t,T\right)$. Therefore, by the elementary independence lemma (e.g., Lemma 2.3.4 in \cite{ShreveFinance}), this means that
\begin{equation*}
\widetilde{c}\left( t,{W}_{t}\right) =\int_{-\infty }^{\widetilde{d%
}_{+}\left({W}_{t},t, T\right) }\frac{1}{\sqrt{2\pi }}e^{-%
\frac{1}{2}y^{2}}dy.
\end{equation*}%
If we apply $\mathbb{D}_{\widetilde{W}_{t}}$ to both sides of the equation then we have by the QCD chain rule (Theorem 3.1 in \cite{A1}) and by the QCD invariance under change of measure that 
\begin{equation}\lbl{dercndexpcttilde}
\begin{split}
\mathbb{D}_{\widetilde{W}_{t}}\widetilde{c}\left( t,{W}_{t}\right) 
&=\frac{1}{\sqrt{2\pi \left( T-t\right) }}e^{-\frac{1}{2}\left( \frac{W_{t}-\int_{t}^{T}\lambda _{s}ds-K}{\sqrt{T-t}}\right) ^{2}}\\
&=p\left(T-t,W_{t}-\int_{t}^{T}\lambda_{s}ds-K\right),
\end{split}
\end{equation}
almost surely, where $p(t,x)$ is the Normal density as in \notnref{papernot}.  This means that, under our assumptions, the digital portfolio is given by 
\begin{equation}\label{DigitalReplicatingPortfolio}
\bsp
\Delta _{t}=e^{-\int_{t}^{T}r(s)ds}a_{t}^{-1}P_{t}^{-1}p\left(T-t,W_{t}-\int_{t}^{T}\lambda_{s}ds-K\right).  
\end{split}
\end{equation}%
Many option prices that are not standard Malliavin differentiable (not in $\sDot$) may be handled similarly.

\subsection{Identifying the chaos expansion of the Brownian indicator}
We now apply the QCD variant of Stroock formula in \thmref{AStr} together with either Theorem 2.1 (i) in \cite{A1} or Theorem 3.1 in \cite{A1} (the QCD fundamental theorem of calculus or chain rule) to identify the integrands of the chaos expansion of $F=\mathbb{I}_{[K,\infty )}\left( W_{T}\right) $. 

First, we make the following simplifying observation about the normal density and its partial derivatives in $x$ given by \eqref{Normaldender}.  
\blm\lbl{Nrmmart}
Assume that the Brownian motion $W$ starts at $x\in\R$.  Then, for any given $y\in\R$ and $n\in\{0,1,2,\ldots\}$, $p_{2}^{(n)}(T-t,W_{t},y)$ is a $W$-martingale for $0\le t<T$ and
\beq\lbl{denbmart}
p_{2}^{(n)}(T-t,W_{t},y)-p_{2}^{(n)}(T,x,y)=\int_{0}^{t}p_{2}^{(n+1)}\lpa T-r, W_{r},y \rpa dW_{r}; \ 0\le t<T,
\eeq
almost surely.  In particular, $\E p_{2}^{(n)}(T-t,W_{t},y)=p_{2}^{(n)}(T,x,y)$ for any $0\le t<T$ and 
$\DWt p_{2}^{(n)}(T-t,W_{t},y)=p_{2}^{(n+1)}\lpa T-t, W_{t},y \rpa$ for all $0\le t<T$ almost surely.
\elm
\bpf
First, we use induction to establish the simple fact that 
\beq\lbl{sf}
\partial_{1}p_{2}^{(n)}\lpa t, x,y\rpa=\frac12\partial^{2}_{22}p_{2}^{(n)}\lpa t, x,y\rpa,\ \forall\ n=0,1,2,\ldots.
\eeq
The assertion is trivially true for $n=0$ since $p(t,x,y)$ is the fundamental solution to the heat equation; in particular,  $\partial_{1}p\lpa t, x,y\rpa=\frac12\partial^{2}_{22}p\lpa t, x,y\rpa$.
Fix an arbitrary $n\in\{0,1,2,\ldots\}$, and assume \eqref{sf} holds for $n$.  We then have 
\beq\lbl{induc}
\bsp
\partial_{1}p_{2}^{(n+1)}\lpa t, x,y\rpa&=\partial_{2}\partial_{1}p_{2}^{(n)}\lpa t, x,y\rpa=\partial_{2}\lbk\frac12\partial^{2}_{22}p_{2}^{(n)}\lpa t, x,y\rpa\rbk\\&=\frac12\partial^{2}_{22}p_{2}^{(n+1)}\lpa t, x,y\rpa,
\end{split}
\eeq
proving \eqref{sf} for every $n\in\{0,1,2,\ldots\}$.

Now, by It\^o's rule and \eqref{sf} we have that for any given $y\in\R$ and $n\in\{0,1,2,\ldots\}$
\beq\lbl{Itodender} 
\bsp
p_{2}^{(n)}(T-t,W_{t},y)-p_{2}^{(n)}(T,x,y)&=-\int_{0}^{t}\partial_{1}p_{2}^{(n)}\lpa T-r, W_{r},y\rpa dr
\\&+\int_{0}^{t}\partial_{2}p_{2}^{(n)}\lpa T-r, W_{r},y\rpa dW_{r}
\\&+\frac12\int_{0}^{t}\partial^{2}_{22}p_{2}^{(n)}\lpa T-r, W_{r},y\rpa dr
\\&=\int_{0}^{t}p_{2}^{(n+1)}\lpa T-r, W_{r},y \rpa dW_{r}.
\end{split}
\eeq
The expectation assertion is trivially obtained by taking expectations on both sides of \eqref{Itodender}, and the $\DW$ assertion follows either by applying $\DWt$ to both sides of \eqref{Itodender} and using Theorem 2.1 (i) in \cite{A1} (the QCD fundamental theorem of calculus) or by applying $\DWt$ to $p_{2}^{(n)}(T-t,W_{t},y)$ and using Theorem 3.1 in \cite{A1} (the QCD chain rule).
\epf
We immediately get the following corollary
\bcr[The chaos expansion of the Brownian indicator]\lbl{chaosidentif}
The Brownian indicator for a Brownian motion $W$ starting at $x\in\R$ has chaos expansion 
$$F:=\mathbb{I}_{[K,\infty )}\left( W_{T}\right)=\sum_{n=0}^{\infty}J_{n}\lpa g_{n}\rpa,$$
where $J_{0}(g_{0})=g_{0}=\E F=\mathbb{P}\left[ W_{T}>K\right]$ and $g_{n}(t_{1},\ldots,t_{n})=p_{2}^{(n-1)}(T,x-K)=\partial_{2}^{n-1}p(T,x-K)$ for all $(t_{1},\ldots,t_{n})\in\left\{ \left( t_{1},t_{2},...,t_{n}\right);0< t_{1}\leq t_{2}\leq ...\leq t_{n}< T\right\}$ and $n\ge1$.
\ecr
\bpf
Recall that by \propref{PropositionIndicator} and by the QCD fundamental theorem of stochastic calculus, Theorem 2.1 in \cite{A1}, we have 
\begin{equation*}
\mathbb{D}_{W_{t}}\mathbb{E}\left[ \mathbb{I}_{[K,\infty )}\left(
W_{T}\right) |\mathscr{F}_{t}\right] =\frac{\exp \left[ \frac{-\left( W_{t}-K\right) ^{2}}{2\left( T-t\right) }\right]}{\sqrt{2\pi \left( T-t\right) }}=p\left(T-t, W_{t}-K\right);\ 0< t<T,
\end{equation*}
almost surely.  Then, by \thmref{AStr} and iterated use of \lemref{Nrmmart} we get
\beq\lbl{coeff}
\bsp
g_{n}(t_{1},\ldots,t_{n})&=\E\lbk\DWtone\E\lbk\ldots\DWtn\E\lbk \mathbb{I}_{[K,\infty )}\left( W_{T}\right)\lgab\sFtn\rbk\cdots\lgab\sFto\rbk\rbk.
\\&=\E p_{2}^{(n-1)}(T-t_{1},W_{t_{1}}-K)=p_{2}^{(n-1)}(T,x-K),
\end{split}
\eeq
for $n\ge1$, and the statement for $n=0$ is trivial, proving our claim.
\epf
Three observations are worth making here
\ben
\rencomrom
\item even though the Brownian indicator $\mathbb{I}_{[K,\infty )}\left( W_{T}\right)$ is not classically Malliavin differentiable (not in $\sDot$), we can easily use the quadratic covariation differentiation theory in \cite{A1} and \thmref{AStr} to obtain its chaos expansion coefficients $g_{n}$ for every $n$.
\item of course, \eqref{coeff} may be rewritten in terms of Hermite polynomials by realizing that if the $n^{th}$ Hermite polynomial is defined by 
\begin{equation*}
H_{n}\left( x\right) =\frac{\left( -1\right) ^{n}}{\sqrt{n!}}e^{\frac{x^{2}}{2%
}}\frac{d^{n}}{dx^{n}}\left( e^{\frac{-x^{2}}{2}}\right);\text{ }n\geq 1,
\end{equation*}%
then \begin{equation*}
p_{2}^{\left( n\right) }\left(t, x\right) =\left( -1\right) ^{n}%
\sqrt{n!}t^{-n/2}p\left(t, x\right) H_{n}\left( 
\frac{x}{\sqrt{t }}\right) ,\text{ }n\geq 1
\end{equation*}%

\item equation \eqref{DWnthder} in \remref{infdiffDW} gives us 
$$\DWtn\E\lbk \mathbb{I}_{[K,\infty )}\left( W_{T}\right)\lgab\sFtn\rbk=p\left(T-t_{n}, W_{t_{n}}-K\right),$$
and \eqref{denbmart} obviates the need for any further dealing with the $\DW$ derivative of the conditional expectations in \eqref {coeff} (since $p_{2}^{(n)}(T-t,W_{t},y)$ is a martingale in $t$ for every $n$); and we were able to apply $\DWt$ to $p\left(T-t, W_{t}-K\right)$ and its spatial derivatives directly using Theorem 3.1 in \cite{A1} (the QCD chain rule).
\een
We now end with a differentiating under the conditional expectation result.
\subsubsection{A conditional QCD chain rule}
In many cases, we have
$$\mathbb{D}_{W_{t}}\mathbb{E}\left[ f\left( W_{T}\right) |\mathscr{F}_{t}\right] =\mathbb{E}\left[ f^{\prime }\left( W_{T}\right) |\mathscr{F}_{t}\right], \forall\  0< t<T. $$
In this subsubsection, we give such a conditional QCD chain rule, whose proof follows from the QCD fundamental theorem of stochastic calculus and the  QCD chain rule given by Theorem 2.1 and Theorem 3.1 in \cite{A1}, respectively.  Now, setting $n=0$ in \eqref{Itodender}  followed by Theorem 2.1 in \cite{A1} (the QCD fundamental theorem of stochastic calculus), and noting that $p_{2}^{(1)}(t,x,y)=-\frac{\partial}{\partial y} p\left(t,x,y\right)$, we obtain that---almost surely---the pair $p\left(T-t,W_{t},y\right)$ and $\DWt p\left(T-t,W_{t},y\right)$ are continuous in $t$ on $[0,T-\epsilon]$ for every $0<\epsilon<T$ and are given by
\beq\lbl{derden}
\bsp
p\left(T-t,W_{t},y\right)-p(T,x,y)&=\int_{0}^{t}-\frac{\partial}{\partial y} p\left(T-r,W_{r},y\right)dW_{r},
\\ \DWt p\left(T-t,W_{t},y\right)&=-\frac{\partial}{\partial y} p\left(T-t,W_{t},y\right)
\\&={\frac {\left(y-W_{t}\right) }{\sqrt {2\pi \, \lpa T-t \rpa^{3} }}}{{\rm e}^{-{\frac { \left(W_{t}-y \right)^{2}}{ 2\lpa T-t \rpa }}}},
\end{split}
\eeq
Another needed ingredient is a standard stochastic Fubini result, which we specialize to our situation and state for the convenience of the reader (see Dol\'ean-Dade \cite{Dol-Da}, Jacod \cite{Jacod} (Th\'eor\`eme 5.44),  and van Neerven et al.~\cite{VanNeer} for more general statements and proofs).
\blm\lbl{sFubini}
Suppose $X:[0,T)\times\R\times\Omega\to\R$ is  $\mathscr{B}([0,T))\times\mathscr{B}(\R)\times \mathscr{F}$-measurable.
If $X_{y}:=\lbr X(s,y,\omega);\ (s,\omega)\in[0,T)\times\Omega\rbr$ is in $\sPtwoslocrTn$ for every $y\in\R$ and every $0< T_{0}<T$; if $\int_{\R}X(s,y,\omega)dy\in\sPtwoslocrTn$ for every $0<T_{0}<T$; and if
\beq\lbl{sFubcnd}
\int_{-\infty }^{\infty}\lab \int_{0}^{t}X(s,y,\omega) dW_{s}\rab dy<\infty; \forall\ 0<t<T, \mbox{ almost surely }
\eeq
then
\beq\lbl{sFubcncl}
\bsp
&\int_{-\infty }^{\infty}\int_{0}^{t}X(s,y,\omega) dW_{s}dy
\\&=\int_{0}^{t}\int_{-\infty }^{\infty}X(s,y,\omega) dydW_{s};\forall\ 0<t<T, \mbox{ almost surely.}
\end{split}
\eeq
 \elm 
 As a corollary, we get
 \bcr\lbl{sFubinicor}
 Suppose $X(s,y,\omega)= f\left( y\right)\mathbb{D}_{W_{s}}p\left(T-s,W_{s},y\right)$, where $f:\R\to\R$ is a continuous function such that $\int_{0}^{t}\lbk\int_{\R}X(s,y,\omega)dy\rbk^{2}ds<\infty$ $\forall\ 0<t<T$ almost surely and $\E\lab f(W_{T})\rab<\infty$.  Then, the conditions of \lemref{sFubini}, including \eqref{sFubcnd}, are satisfied and 
 \beq\lbl{sFubcncl2}
 \bsp
&\int_{-\infty }^{\infty}\int_{0}^{t}f\left( y\right)\mathbb{D}_{W_{s}}p\left(T-s,W_{s},y\right) dW_{s}dy
\\&=\int_{0}^{t}\int_{-\infty }^{\infty}f\left( y\right)\mathbb{D}_{W_{s}}p\left(T-s,W_{s},y\right) dydW_{s};\forall\ 0<t<T, \mbox{ almost surely.}
\end{split}
\eeq
 \ecr
\bpf  We only need to verify \eqref{sFubcnd}.  The $L^{1}$ condition $\E\lab f(W_{T})\rab<\infty$, the Markov property, \eqref{derden}, and the fact that $p\left(t,x,y\right)>0$ for all $t$, $x$, and $y$ mean that 
\beq\lbl{verifyFub}
\bsp
\infty&>\E\lbk\lab f\lpa W_{T} \rpa\rab|\sFt\rbk+\E\lab f\lpa W_{T} \rpa\rab=\int_{-\infty }^{\infty }\lab f\left( y\right)\rab p\left(T-t,W_{t},y\right) dy\\&+\int_{-\infty }^{\infty }\lab f\left( y\right)\rab p\left(T,x,y\right) dy
\\&\ge\int_{-\infty }^{\infty }\lab f\left( y\right)\rab\lab\int_{0}^{t} \DWs p\left(T-s,W_{s},y\right)dW_{s}\rab dy
\\&=\int_{-\infty }^{\infty} \lab\int_{0}^{t}f\left( y\right) \DWs p\left(T-s,W_{s},y\right)dW_{s}\rab dy
\end{split}
\eeq
for all $0<t<T$ almost surely, verifying \eqref{sFubcnd}.
\epf
\blm[Conditional QCD chain rule]\label{ConditionalChainRuleLemma}
Assume that $f\in C^{1}\left(\mathbb{R};\mathbb{R}\right) $ and that
\beq\lbl{cnd117}
\bsp
&\int_{0}^{t}\lbk\int_{\R}f\left( y\right)\mathbb{D}_{W_{s}}p\left(T-s,W_{s},y\right)dy\rbk^{2}ds<\infty;\ 0<t<T \mbox{ almost surely,}\\
&\mathbb{E}\left\vert f'\left( W_{T}\right) \right\vert <\infty\mbox{ and }\E\lab f(W_{T})\rab<\infty.
\end{split}
\eeq
Then, almost surely,
\begin{equation*}
\mathbb{D}_{W_{t}}\mathbb{E}\left[ f\left( W_{T}\right) |\mathscr{F}_{t}\right] =\mathbb{E}\left[ f^{\prime }\left( W_{T}\right) |\mathscr{F}_{t}\right], \forall\  0< t<T. 
\end{equation*}
\elm
\bpf
By the Markov property, integration by parts, equation \eqref{derden}, Theorem 2.1 (i) in \cite{A1}, \coref{sFubinicor}, \eqref{derden} again, and the Markov property again, we obtain  
\beqs
\begin{split}
&\mathbb{E}\left[ f^{\prime }\left(W_{T}\right) |\mathscr{F}_{t}\right]=\int_{-\infty }^{\infty }f'\left( y\right)p\left(T-t,W_{t},y\right) dy=-\int_{-\infty }^{\infty }f\left( y\right) \frac{\partial }{\partial y}p\left(T-t,W_{t},y\right) dy
\\&=\int_{-\infty }^{\infty }f\left( y\right) \mathbb{D}_{W_{t}}p\left(T-t,W_{t},y\right) dy=\mathbb{D}_{W_{t}}\int_{0}^{t}\int_{-\infty }^{\infty }f\left( y\right)\mathbb{D}_{W_{s}}p\left(T-s,W_{s},y\right) dydW_{s}
\\&=\mathbb{D}_{W_{t}}\int_{-\infty }^{\infty }f\left( y\right) \int_{0}^{t}\mathbb{D}_{W_{s}}p\left(T-s,W_{s},y\right) dW_{s}dy \\ 
&=\mathbb{D}_{W_{t}}\lbk\int_{-\infty }^{\infty }f\left( y\right) p\left(T-t,W_{t},y\right) dy-\int_{-\infty }^{\infty }f\left( y\right) p\left(T,x,y\right) dy\rbk\\&=
\mathbb{D}_{W_{t}}\mathbb{E}\left[ f\left( W_{T}\right) |\mathscr{F}_{t}\right] 
\end{split}
\eeqs
for every  $ 0<t<T $ almost surely.  This completes the proof of our assertion.
\epf
\brm
It is also easy to apply \thmref{AStr} (b), \eqref{dercndexpcttilde}, and \lemref{ConditionalChainRuleLemma} to get the chaos expansion coefficients in $\mathbb{I}_{[K,\infty )}=\sum_{n=0}^{\infty }\widetilde{J}_{n}\left( g_{n}\right),$ where $\widetilde{J_{n}}$ denotes the $n$-fold iterated Wiener integral with respect to $\widetilde{W}$ over the set $S_{n}$.  We leave this to the interested reader.
\erm
\appendix\section{Notation, definitions, and a brief review of Girsanov's theorem and the quadratic covariation of processes}\lbl{appA}
\bnt\lbl{papernot}
For typesetting and aesthetic reasons, we alternate freely between $X_{t}$ and $X(t)$ to denote any stochastic process $X$ evaluated at $t$.  We denote by $p(t,x,y)$ the Normal density ${\exp \left[ \frac{-(x-y)^{2}}{2t}\right]}/{\sqrt{2\pi t }}$, and we define $p(t,x)=p(t,x,0)$.  Also, for  a function $f$ of $d$ variables $x_{1},\ldots,x_{d}$ we denote the $n$th order derivative in the $k$th variable by $\partial^{n}_{k}f(x_{1},\ldots,x_{d})$.  I.e.,
\beq\lbl{partialdef}
\partial^{n}_{k}f(x_{1},\ldots,x_{d})=\frac{\partial^{n}}{\partial x_{k}^{n}}f(x_{1},\ldots,x_{d});\ k=1,\ldots,d,\ n=0,1,2,\ldots
\eeq
with $\partial^{0}_{k}f(x_{1},\ldots,x_{d})=f(x_{1},\ldots,x_{d})$ for any $k=1,\ldots,d$.  Finally, we use
\beq\lbl{Normaldender}
\bsp
&p_{1}^{(n)}(t,x,y):=\partial^{n}_{1}p(t,x,y)=\frac{\partial^{n}}{\partial t^{n}}\frac{\exp \left[ \frac{-(x-y)^{2}}{2t}\right]}{\sqrt{2\pi t }}; \ n=0,1,2,\ldots
\\& p_{2}^{(n)}(t,x,y):=\partial^{n}_{2}p(t,x,y)=\frac{\partial^{n}}{\partial x^{n}}\frac{\exp \left[ \frac{-(x-y)^{2}}{2t}\right]}{\sqrt{2\pi t }}; \ n=0,1,2,\ldots
\end{split}
\eeq
\ent
\begin{defn}[Integrand Classes]\label{intgrnds}
Let $\sPtworT$ be the class of processes $X:[0,T]\times \Omega \rightarrow \mathbb{R}$ such that
\begin{enumerate}\renewcommand{\labelenumi}{$(\arabic{enumi})$}
\item $X$ is measurable: $\left( t,\omega \right) \rightarrow X\left( t,\omega \right) $ is $\mathscr{B}([0,T])\times 
\mathscr{F}$ measurable, where $\mathscr{B}([0,T])$ denotes the Borel $\sigma$-algebra on $[0,T],$
\item $X$ is $\lbr\sFt\rbr$-adapted: $X_{t}\in\sFt$, for every $t\in[0,T],$
\item $\ds\mathbb{E}\left[ \int_{0}^{T}X^{2}\left( s,\omega \right) ds\right]<\infty.$ 
\end{enumerate}
The class $\sPtwosrT\subset\sPtworT$ is obtained from $\sPtworT$ by leaving condition $(3)$ unchanged and replacing the measurability and adaptability requirements in $(1)$ and $(2)$ by the stronger requirement of progressive measurability with respect to the filtration $\lbr\sFt\rbr$$:$
$$\left( t,\omega \right) \rightarrow X\left( t,\omega \right)\mbox{ is } \mathscr{B}([0,t])\times\sFt\mbox{ measurable for each }t\in[0,T].$$
  The classes $\sPtwolocrT$ and $\sPtwoslocrT$ are obtained from $\sPtworT$ and $\sPtwosrT$, respectively, by replacing condition $(3)$ with the weaker condition
\beq\lbl{locint}
\P\lbk\int_{0}^{T}X^{2}\left( s,\omega \right) ds<\infty\rbk=1.
\eeq  
 \end{defn} 
 \bdf[Almost indistinguishability and indistinguishability]\lbl{aindisdef} Suppose $X,Y:\Omega\times[0,T]\to\R$ are two stochastic processes on a probability space $\OFP$ such that 
\beq\lbl{aindis}
\bsp
Y_{t}\lpa\omega\rpa=X_{t}\lpa\omega\rpa;\quad t\in[0,T]\setminus Z(\omega),\ \omega\in\Omega^{*}.
\end{split}
\eeq
holds for some subset $\Omega^{*}\subset\Omega$, with $\P\lpa \Omega^{*}\rpa=1$ for a collection of zero Lebesgue-measure random sets $\lbr Z(\omega); \omega\in\Omega^{*}\rbr$.  Then, we say that $X$ and $Y$ are almost indistinguishable versions of each other.  $X$ and $Y$ are indistinguishable if $Z(\omega)=\phi$ (the empty set) for each $\omega\in\Omega^{*},$ and we write $X=Y$ a.s.  If $Y_{t}(\omega)$ is defined only for $(\omega,t)$ where $t\in[0,T]\setminus Z(\omega)$ and $\omega\in\Omega^{*}$ and if the stochastic process $X:\Omega\times[0,T]\to\R$ satisfies \eqref{aindis}, then $X$ is said to be an almost indistinguishable extension of $Y$.  Any such extension $X$ of $Y$ is denoted by $Y^{\sc{aie}}$.
\edf
\brm
It is clear that two stochastic processes $X,Y:\Omega\times[0,T]\to\R$ on $\OFP$ are almost indistinguishable versions of one another iff 
\beq\lbl{aiinteg}
\int_{0}^{T}\lab X(t)-Y(t)\rab ^{2}dt=0,\text{ a.s. }\P.
\eeq
It is also obvious that two almost indistinguishable extensions of $Y$, $Y^{\sc{aie}}_{1}$ and $Y^{\sc{aie}}_{2}$ are almost indistinguishable versions of one another and hence satisfy \eqref{aiinteg}.  
\erm

We assume the same setup as the classical Girsanov change of measure theorem, which we now combine with the subsequently discovered Novikov sufficient condition (see e.g., \cite{KaratzasShreve,ShreveFinance}). 
\begin{thm}[Girsanov 1960 and Novikov 1972]\label{Girsanov}
Let $\lambda =\left\{ \lambda _{t},\mathscr{F}_{t}:t\in \left[ 0,T\right] \right\}\in\sPtwoslocrT.$ Define 
\begin{equation*}
Z_{t}=\exp \left\{ -\int_{0}^{t}\lambda \left( u\right) dW\left( u\right) -%
\frac{1}{2}\int_{0}^{t}\lambda ^{2}\left( u\right) du\right\} ,
\end{equation*}%
\begin{equation*}
\widetilde{W}_{t}=W_{t}+\int_{0}^{t}\lambda \left( u\right) du; \ 0\le t\le T
\end{equation*}%
and assume that either one of the two following conditions hold
\begin{equation}\label{GirsanovCondition}
(i)\ \mathbb{E}\int_{0}^{T}\lambda _{u}^{2}Z_{u}^{2}du<\infty \mbox{ or }(ii)\ \E\lbk\exp\lpa\frac12\int_{0}^{T}\lambda^{2} _{s}ds\rpa\rbk<\infty.
\end{equation}%
Then $ Z=\lbr Z_{t},\sFt;0\le t\le T\rbr$ is a martingale, $\mathbb{E}\left[Z_{t}\right] =1$ for $0\le t\le T$, and if $\widetilde{\mathbb{P}}$ is defined by the recipe 
\begin{equation*}
\frac{d\widetilde{\mathbb{P}}}{d\mathbb{P}}=Z_{T},
\end{equation*}%
then $\Pt$ is a probability measure on $\sFT$ and the process $\widetilde{W}=\lbr \widetilde{W}_{t},\sFt;0\le t\le T\rbr$ is a Brownian motion on the probability space $\OFTFtPt$.
\end{thm}

We denote by $\Et$ the expectation taken with respect to $\Pt$.  We use the following standard result regarding the behavior of conditional expectations under change of measure (see page 193 of \cite{KaratzasShreve})

\begin{lem}[Bayes Rule]\label{CondExpCOM}
If $0\leq s\leq t\leq T$ and $F$ is an $\mathscr{F}_{t}$-measurable random variable such that $\Et\lab F\rab<\infty$. If $Z$ is a martingale, then 
\begin{equation*}
\widetilde{\mathbb{E}}\left[ F|%
\mathscr{F}_{s}\right] =\frac{1}{Z_{s}}\mathbb{E}\left[ FZ_{t}|\mathscr{F}_{s}\right]\mbox{ a.s. }\P\mbox{ and }\Pt. 
\end{equation*}
\end{lem}
We now recall the definition of the covariation process of two processes.  We denote by $\stackrel{\P}{\to}$ convergence in probability under the probability measure $\P$.

\bdf\lbl{covXY}
Two real-valued processes $X,Y$ on a probability space $\OFP$ have finite quadratic covariation iff there exists a finite process  $\lqv X,Y\rqv$ such that for every $t>0$ and every sequence $\lbr \T_{n}\rbr$ of partitions of $[0,t]$---$\T_{n}=\lbr t_{0},t_{1},\ldots,t_{n}\rbr$ with $0=t_{0}<t_{1}<\cdots<t_{n}=t$---such that the mesh limit $\lim_{n\to\infty}\lab\T_{n}\rab=0$
\beq\lbl{qcovdef}
V^{2}_{t}\lpa X,Y,\T_{n}\rpa:=\sum_{k=1}^{n}\lpa X_{t_{k}}-X_{t_{k-1}}\rpa\lpa Y_{t_{k}}-Y_{t_{k-1}}\rpa\stackrel{\P}{\to}\lqv X,Y\rqv_{t}\mbox{ as }n\to\infty.
\eeq
The process $\lqv X,Y\rqv$ is called the quadratic covariation of $X$ and $Y$.  The process $\lqv X,X\rqv$ (the case $X\equiv Y$) is called the quadratic variation of $X$.  
 \edf
 When we want to emphasize the role of $\P$ in the definition of $\lqv X,Y\rqv$, we write $\lqv X,Y\rqv^{\P}$.
It is then a simple matter to see the following invariance-under-equivalent-change-of-measure property of the process $\lqv\cdot,\cdot\rqv$.
\blm\lbl{invcomqcp}
Let $T>0$ be fixed but arbitrary.  Suppose $X$ and $Y$ are two real-valued adapted processes defined on the interval $[0,T]$ and on the probability space $\OFFtP$$;$ and suppose that a probability measure $\Pt$ is defined on $\sFT$ and is equivalent to the restriction of $\P$ to $\sFT$.  If either one of $\lqv X,Y\rqv^{\P}$ or $\lqv X,Y\rqv^{\Pt}$ is finite on $[0,T]$, then so is the other and
\beq\lbl{invcomqcpeq}
\lqv X,Y\rqv^{\P}_{t}=\lqv X,Y\rqv^{\Pt}_{t}; \ \mbox{ a.s. }\P \mbox{ and } \Pt,\ \forall\ 0< t\le T.
\eeq
I.e., they are modifications of one another under both $\P$ and $\Pt$.   In particular, if $\Pt$ is the Girsanov probability measure in \thmref{Girsanov} with $\lambda$ satisfying \eqref{GirsanovCondition}((i) or (ii)) and if $X$ and $Y$ are continuous semimartingales, then $\lqv X,Y\rqv^{\P}$ and $\lqv X,Y\rqv^{\Pt}$ are indistinguishable under both $\P$ and $\Pt$$:$ $$\P\lbk\lqv X,Y\rqv^{\P}_{t}=\lqv X,Y\rqv^{\Pt}_{t};0\le t\le T\rbk=1=\Pt\lbk\lqv X,Y\rqv^{\P}_{t}=\lqv X,Y\rqv^{\Pt}_{t};0\le t\le T\rbk.$$
\elm
\bpf
Let $Z_{n,t}:=V^{2}_{t}\lpa X,Y,\T_{n}\rpa$; then by \defnref{covXY}, elementary measure theory, and the equivalence of $\P$ and $\Pt$ we have 
\beqs
\bsp
&Z_{t}=\lqv X,Y\rqv^{\P}_{t} \mbox{ is finite }\iff Z_{n,t}\stackrel{\P}{\to}Z_{t}\mbox{ as }n\to\infty\\
&\iff\mbox{for every subsequence }\lbr Z_{n_{k},t}\rbr_{k}\mbox{ there is a further subsequence }\lbr Z_{n_{k_{l}},t}\rbr_{l}\\
&\mbox{such that }\P\lbk \lim_{l\to\infty}  Z_{n_{k_{l}},t}=Z_{t}\rbk=1=\Pt\lbk \lim_{l\to\infty}  Z_{n_{k_{l}},t}=Z_{t}\rbk
\\&\iff Z_{n,t}\stackrel{\Pt}{\to}Z_{t}\mbox{ as }n\to\infty.
\end{split}
\eeqs
for every fixed $0< t\le T$ and \eqref{invcomqcpeq} follows.  If $X$ and $Y$ are continuous semimartingales, $\lambda$ satisfies \eqref{GirsanovCondition} ((i) or (ii)), and $\Pt$ is Girsanov's probability measure given in \thmref{Girsanov}; then both $\lqv X,Y\rqv^{\P}$ and $\lqv X,Y\rqv^{\Pt}$ are modifications of one another and are both almost surely ($\P$ and $\Pt$) continuous and hence indistinguishable under both $\P$ and $\Pt$. 
 \epf
\section{Briefly on chaos expansion and standard Malliavin's derivative versions of Clark-Ocone Formulas}\lbl{appB}
 Let $F\in L^{2}\left(\Omega,\mathscr{F}_{T}, \mathbb{P}\right)$ be an $L^{2}$ and $\mathscr{F}_{T}$-measurable random variable, with chaos
expansion 
\begin{equation*}
F=\sum_{n=0}^{\infty }I_{n}\left( \hat{f}_{n}\right),
\end{equation*}%
where $I_{n}$ is the $n$-fold iterated It\^o-Wiener integral over $[0,T]^{n}$ 
\begin{equation*}
I_{n}\left( \hat{f}_{n}\right) =\int_{\left[ 0,T\right] ^{n}}\hat{f}_{n}\left( t_{1},...,t_{n}\right) dW_{t_{1}}dW_{t_{2}}...dW_{t_{n-1}}dW_{t_{n}}\end{equation*}
and $\hat{f}_{n}\in \hat{L}^{2}\left([0,T] ^{n}\right)$, where $\hat{L}^{2}\left([0,T]^{n}\right) $ is the space of symmetric Borel
deterministic square integrable functions. Then, we have the isometry 
\beq\lbl{L2Fiso}
\lnrm F\rnrm^{2}_{L^{2}(\Omega,\P)}=\sum_{n=0}^{\infty}n!\lnrm \hat{f}_{n}\rnrm^{2}_{L^{2}([0,T]^{n})}=\sum_{n=0}^{\infty}\E\lbk I_{n}\lpa\hat{f}\rpa\rbk^{2}
\eeq
We say that $F\in\sDot$ if 
\begin{equation}\lbl{sDotnrm}
\left\Vert F\right\Vert _{\sDot}^{2}:=\sum_{n=1}^{\infty }nn!\left\Vert 
\hat{f}_{n}\right\Vert _{L^{2}\left([0,T]^{n}\right)
}^{2}<\infty .
\end{equation}
For $F\in\sDot$, we define the Malliavin derivative $D_{t}F$ of $F$ at time $t$ as the expansion
$$D_{t}F=\sum_{n=1}^{\infty }nI_{n-1}\left( \hat{f}_{n}(\cdot,t)\right);\quad t\in[0,T],$$ 
where $I_{n-1}\left( \hat{f}_{n}(\cdot,t)\right)$ is the $(n-1)$ fold iterated integral of $\hat{f}(t_{1},\ldots,t_{n-1},t)$ with respect to the first $n-1$ variables $t_{1},\ldots,t_{n-1}$ and $t_{n}=t$ left as a parameter.  Observe that 
\beq\lbl{l2vsd12}
\bsp
\lnrm F\rnrm_{\sDot}^{2}&=\sum_{n=1}^{\infty}\int_{0}^{T}n^{2}(n-1)!\lnrm  \hat{f}_{n}(\cdot,t)\rnrm _{L^{2}\lpa[0,T]^{n-1}\rpa}^{2}dt\\
&=\int_{0}^{T}\E\sum_{n=1}^{\infty}n^{2}\lbk I_{n-1}\left( \hat{f}_{n}(\cdot,t)\right)\rbk^{2}dt=\int_{0}^{T}\E\lpa D_{t}F\rpa^{2}dt
\\&=\lnrm D_{\cdot}F\rnrm^{2}_{L^{2}\lpa\Omega\times[0,T],\P\times\lambda\rpa},
\end{split}
\eeq
where we used the fact that $\E\lbk I_{n}(g)\rbk^{2}=n!\lnrm g\rnrm^{2}_{L^{2}\lpa[0,T]^{n}\rpa}$ along with the isometry in \eqref{L2Fiso} and where $\lambda$ is Lebesgue's measure on $[0,T]$.
\begin{thm}[The Standard Malliavin derivative Clark-Ocone formula] \lbl{sCO}Let $F\in D_{1,2}$ be $\mathscr{F}_{T}$-measurable. Then 
\begin{equation*}
F=\mathbb{E}\left[ F\right] +\int_{0}^{T}\mathbb{E}\left[ D_{t}F|%
\mathscr{F}%
_{t}\right] dW_{t}; \mbox{ a.s. }\P.
\end{equation*}%
\end{thm}
The standard Clark-Ocone formula under change of measure (COM), using Malliavin's derivative, was introduced in \cite{OcKaCOM}. For the sake
of comparison we include it in the next theorem (see \cite{Levy} page 46).

\begin{thm}[The standard Malliavin derivative Clark-Ocone formula under COM] 
\label{COCOM}%
Suppose $F$ is $%
\mathscr{F}%
_{T}$ measurable, $F\in\sDot$, and that 
\begin{equation}\lbl{COCOMcnd1}
\widetilde{\mathbb{E}}\left[ \left\vert F\right\vert \right] <\infty ,
\end{equation}%
\begin{equation}\lbl{COCOMcnd2}
\widetilde{\mathbb{E}}\left[ \int_{0}^{T}\left\vert D_{t}F\right\vert ^{2}dt%
\right] <\infty ,\text{ }
\end{equation}%
\begin{equation}\lbl{COCOMcnd3}
\widetilde{\mathbb{E}}\left[ \left\vert F\right\vert \int_{0}^{T}\left(
\int_{0}^{T}D_{t}\lambda _{s}dW_{s}+\int_{0}^{T}\left( D_{t}\lambda
_{s}\right) \lambda _{s}ds\right) ^{2}dt\right] <\infty ,
\end{equation}%
where $D_{t}$ is the Malliavin derivative. Then, 
\begin{equation}\lbl{COCOMcncl}
F=\widetilde{\mathbb{E}}\left[ F\right] +\int_{0}^{T}\widetilde{\mathbb{E}}%
\left[\left.\left( D_{t}F-F\int_{t}^{T}D_{t}\lambda _{s}d\widetilde{W}_{s}\right)\right|%
\mathscr{F}_{t}\right] d\widetilde{W}_{t}.
\end{equation}
\end{thm}
\section{Frequent acronyms and notations key}\lbl{B}
\begin{enumerate}\renewcommand{\labelenumi}{\Roman{enumi}.}
\item {\textbf{Acronyms}}\vspace{2mm}
\begin{enumerate}\renewcommand{\labelenumii}{(\arabic{enumii})}
\item BM: Brownian motion,
\item QCD: quadratic covariation derivative (see Allouba's article \cite{A1}). 
\item RCLL (or cadlag): right continuous with left limits.
\end{enumerate}
\vspace{2.5mm}
\item {\textbf{Notations}}\vspace{2mm}
\begin{enumerate}\renewcommand{\labelenumii}{(\arabic{enumii})}
\item $\DW$:  the QCD process,
\item $\DWt$: the QCD at time $t$,
\item $\Et$: the expectation taken with respect to Girsanov's changed probability measure $\Pt$,
\item $\partial^{n}_{k}$, $p_{1}^{(n)}(t,x,y)$, $p_{2}^{(n)}(t,x,y)$: see \notnref{papernot},
\item $\lqv X, Y\rqv$:  the quadratic covariation processs of the processes $X$  and $Y$ (\defnref{covXY}),
\item $\sPtworT$, $\sPtwosrT$, $\sPtwolocrT$, and $\sPtwoslocrT$:  standard classes of integrands with respect to the BM $W$ on the interval $[0,T]$ (\defnref{intgrnds}),
\item $\S_{n}$: the simplex $\S_{n}=\left\{ \left( t_{1},t_{2},...,t_{n}\right);0\leq t_{1}\leq t_{2}\leq ...\leq t_{n}\leq T\right\}\subset[0,T]^{n}$.
\item $\sDot$:  Malliavin's standard space of differentiable random variables whose Malliavin derivative is in $L^{2}$ (see \eqref{sDotnrm} and \eqref{l2vsd12} for precise statements). 
\end{enumerate}
\end{enumerate}

\end{document}